\documentclass[12pt]{article}
\usepackage{amssymb,latexsym}

\def\Wrfi{W_{r,\underline{f}}}
\def\Wrxi{W_{r,\underline{x}}}
\def\Wdxi{W_{d,\underline{x}}}
\def\Wto{W_{2,1}}
\def\Wrxit{\overline{W}_{r,\underline{x}}}
\def\Wrfit{\overline{W}_{r,\underline{f}}}
\def\sympol{R[x_1,\ldots,x_r]^{S_r}}
\def\sympold{R[x_1,\ldots,x_d]^{S_d}}
\def\deg{\mathop{\rm deg}\nolimits}
\def\dim{\mathop{\rm dim}\nolimits}
\def\sdim{\mathop{\rm sdim}\nolimits}
\def\image{\mathop{\rm image}\nolimits}
\def\span{\mathop{\rm span}\nolimits}
\def\End{\mathop{\rm End}\nolimits}
\def\tr{\mathop{\rm tr}\nolimits}
\def\str{\mathop{\rm str}\nolimits}
\def\Hom{\mathop{\rm Hom}\nolimits}
\def\gr{\mathop{\rm gr}\nolimits}
\def\ad{\mathop{\rm ad}\nolimits}
\def\id{\mathop{\rm id}\nolimits}
\def\Gal{\mathop{\rm Gal}\nolimits}
\def\MO{\mathop{\rm M_0}\nolimits}
\def\hatA{\hat{A}}
\def\hatAtan{\hat{A}_{t}}
\def\Ztan{Z_{t}}
\def\hatAc{\hat{A}^c}
\def\Vkan{{\cal V}^{can}}
\def\K{F}
\def\OverGam{\overline{\Gamma}}
\def\VI{\hat{\Phi}}
\def\V1{\Upsilon}
\def\VS{\Phi}
\def\WS{\Psi}

\def\vs{vector space}
\def\overtilde{\overline}
\def\fknots{{\cal K}}
\def\realize{\psi}
\def\forget{\rho}

\def\L{\mathfrak{g}}
\def\Lnr{\L_n^{\oplus r}}
\def\Wglnxt{W_{\gl_n,X,\tau}}
\def\edge[#1,#2]{]#1;#2[}
\def\gl{\mathfrak{gl}}
\def\so{\mathfrak{so}}

\def\spli{\mathfrak{sl}}

\def\osp{\mathfrak{osp}}
\def\emptyseq{\emptyset}

\newtheorem{lemma}{Lemma}[section]
\newtheorem{prop}{Proposition}[section]
\newtheorem{remark}{Remark}[section]
\newtheorem{theorem}{Theorem}[section]
\newtheorem{coro}{Corollary}[section]
\newtheorem{defi}{Definition}[section]

\newcommand{\N}{{\ensuremath{\mathbb N}}}
\newcommand{\Z}{{\ensuremath{\mathbb Z}}}
\newcommand{\Q}{{\ensuremath{\mathbb Q}}}
\newcommand{\R}{{\ensuremath{\mathbb R}}}

\def\DottedCircle{
\bezier{4}(0.966,-0.259)(1.04,0)(0.966,0.259)
\bezier{4}(0.966,0.259)(0.897,0.518)(0.707,0.707)
\bezier{4}(0.707,0.707)(0.518,0.897)(0.259,0.966)
\bezier{4}(0.259,0.966)(0,1.04)(-0.259,0.966)
\bezier{4}(-0.259,0.966)(-0.518,0.897)(-0.707,0.707)
\bezier{4}(-0.707,0.707)(-0.897,0.518)(-0.966,0.259)
\bezier{4}(-0.966,0.259)(-1.04,0)(-0.966,-0.259)
\bezier{4}(-0.966,-0.259)(-0.897,-0.518)(-0.707,-0.707)
\bezier{4}(-0.707,-0.707)(-0.518,-0.897)(-0.259,-0.966)
\bezier{4}(-0.259,-0.966)(0,-1.04)(0.259,-0.966)
\bezier{4}(0.259,-0.966)(0.518,-0.897)(0.707,-0.707)
\bezier{4}(0.707,-0.707)(0.897,-0.518)(0.966,-0.259)
}
\def\FullCircle{
\thicklines
\put(0,0){\circle{2}}
}
\def\Endpoint[#1]{
\ifcase#1
\put(1,0){\circle*{0.15}}
\or\put(0.866,0.5){\circle*{0.15}}
\or\put(0.5,0.866){\circle*{0.15}}
\or\put(0,1){\circle*{0.15}}
\or\put(-0.5,0.866){\circle*{0.15}}
\or\put(-0.866,0.5){\circle*{0.15}}
\or\put(-1,0){\circle*{0.15}}
\or\put(-0.866,-0.5){\circle*{0.15}}
\or\put(-0.5,-0.866){\circle*{0.15}}
\or\put(0,-1){\circle*{0.15}}
\or\put(0.5,-0.866){\circle*{0.15}}
\or\put(0.866,-0.5){\circle*{0.15}}
\fi}
\def\Arc[#1]{
\thicklines                     
\ifcase#1
\bezier{25}(0.966,-0.259)(1.04,0)(0.966,0.259)
\or
\bezier{25}(0.966,0.259)(0.897,0.518)(0.707,0.707)
\or
\bezier{25}(0.707,0.707)(0.518,0.897)(0.259,0.966)
\or
\bezier{25}(0.259,0.966)(0,1.04)(-0.259,0.966)
\or
\bezier{25}(-0.259,0.966)(-0.518,0.897)(-0.707,0.707)
\or
\bezier{25}(-0.707,0.707)(-0.897,0.518)(-0.966,0.259)
\or
\bezier{25}(-0.966,0.259)(-1.04,0)(-0.966,-0.259)
\or
\bezier{25}(-0.966,-0.259)(-0.897,-0.518)(-0.707,-0.707)
\or
\bezier{25}(-0.707,-0.707)(-0.518,-0.897)(-0.259,-0.966)
\or
\bezier{25}(-0.259,-0.966)(0,-1.04)(0.259,-0.966)
\or
\bezier{25}(0.259,-0.966)(0.518,-0.897)(0.707,-0.707)
\or
\bezier{25}(0.707,-0.707)(0.897,-0.518)(0.966,-0.259)
\fi}
\def\DottedArc[#1]{
\ifcase#1
\bezier{4}(0.966,-0.259)(1.04,0)(0.966,0.259)
\or
\bezier{4}(0.966,0.259)(0.897,0.518)(0.707,0.707)
\or
\bezier{4}(0.707,0.707)(0.518,0.897)(0.259,0.966)
\or
\bezier{4}(0.259,0.966)(0,1.04)(-0.259,0.966)
\or
\bezier{4}(-0.259,0.966)(-0.518,0.897)(-0.707,0.707)
\or
\bezier{4}(-0.707,0.707)(-0.897,0.518)(-0.966,0.259)
\or
\bezier{4}(-0.966,0.259)(-1.04,0)(-0.966,-0.259)
\or
\bezier{4}(-0.966,-0.259)(-0.897,-0.518)(-0.707,-0.707)
\or
\bezier{4}(-0.707,-0.707)(-0.518,-0.897)(-0.259,-0.966)
\or
\bezier{4}(-0.259,-0.966)(0,-1.04)(0.259,-0.966)
\or
\bezier{4}(0.259,-0.966)(0.518,-0.897)(0.707,-0.707)
\or
\bezier{4}(0.707,-0.707)(0.897,-0.518)(0.966,-0.259)
\fi}
\def\Chord[#1,#2]{
\thinlines
\ifnum#1>#2\Chord[#2,#1]
\else\ifnum#1<#2
\ifcase#1
\ifcase#2
\or\qbezier(1,0)(0.516,0.138)(0.866,0.5)
\or\qbezier(1,0)(0.45,0.26)(0.5,0.866)
\or\qbezier(1,0)(0.327,0.327)(0,1)
\or\qbezier(1,0)(0.179,0.311)(-0.5,0.866)
\or\qbezier(1,0)(0.0536,0.2)(-0.866,0.5)
\or\put(1, 0){\line(-2, 0){2}}
\or\qbezier(1,0)(0.0536,-0.2)(-0.866,-0.5)
\or\qbezier(1,0)(0.179,-0.311)(-0.5,-0.866)
\or\qbezier(1,0)(0.327,-0.327)(0,-1)
\or\qbezier(1,0)(0.45,-0.26)(0.5,-0.866)
\or\qbezier(1,0)(0.516,-0.138)(0.866,-0.5)
\fi
\or\ifcase#2\or
\or\qbezier(0.866,0.5)(0.378,0.378)(0.5,0.866)
\or\qbezier(0.866,0.5)(0.26,0.45)(0,1)
\or\qbezier(0.866,0.5)(0.12,0.446)(-0.5,0.866)
\or\qbezier(0.866,0.5)(0,0.359)(-0.866,0.5)
\or\qbezier(0.866,0.5)(-0.0536,0.2)(-1,0)
\or\put(0.866, 0.5){\line(-5, -3){1.73}}
\or\qbezier(0.866,0.5)(0.146,-0.146)(-0.5,-0.866)
\or\qbezier(0.866,0.5)(0.311,-0.179)(0,-1)
\or\qbezier(0.866,0.5)(0.446,-0.12)(0.5,-0.866)
\or\qbezier(0.866,0.5)(0.52,0)(0.866,-0.5)
\fi
\or\ifcase#2\or\or
\or\qbezier(0.5,0.866)(0.138,0.516)(0,1)
\or\qbezier(0.5,0.866)(0,0.52)(-0.5,0.866)
\or\qbezier(0.5,0.866)(-0.12,0.446)(-0.866,0.5)
\or\qbezier(0.5,0.866)(-0.179,0.311)(-1,0)
\or\qbezier(0.5,0.866)(-0.146,0.146)(-0.866,-0.5)
\or\put(0.5, 0.866){\line(-3, -5){1}}
\or\qbezier(0.5,0.866)(0.2,-0.0536)(0,-1)
\or\qbezier(0.5,0.866)(0.359,0)(0.5,-0.866)
\or\qbezier(0.5,0.866)(0.446,0.12)(0.866,-0.5)
\fi
\or\ifcase#2\or\or\or
\or\qbezier(0,1.)(-0.138,0.516)(-0.5,0.866)
\or\qbezier(0,1.)(-0.26,0.45)(-0.866,0.5)
\or\qbezier(0,1.)(-0.327,0.327)(-1,0)
\or\qbezier(0,1.)(-0.311,0.179)(-0.866,-0.5)
\or\qbezier(0,1.)(-0.2,0.0536)(-0.5,-0.866)
\or\put(0, 1){\line(0, -2){2}}
\or\qbezier(0,1.)(0.2,0.0536)(0.5,-0.866)
\or\qbezier(0,1.)(0.311,0.179)(0.866,-0.5)
\fi
\or\ifcase#2\or\or\or\or
\or\qbezier(-0.5,0.866)(-0.378,0.378)(-0.866,0.5)
\or\qbezier(-0.5,0.866)(-0.45,0.26)(-1,0)
\or\qbezier(-0.5,0.866)(-0.446,0.12)(-0.866,-0.5)
\or\qbezier(-0.5,0.866)(-0.359,0)(-0.5,-0.866)
\or\qbezier(-0.5,0.866)(-0.2,-0.0536)(0,-1)
\or\put(-0.5, 0.866){\line(3, -5){1}}
\or\qbezier(-0.5,0.866)(0.146,0.146)(0.866,-0.5)
\fi
\or\ifcase#2\or\or\or\or\or
\or\qbezier(-0.866,0.5)(-0.516,0.138)(-1,0)
\or\qbezier(-0.866,0.5)(-0.52,0)(-0.866,-0.5)
\or\qbezier(-0.866,0.5)(-0.446,-0.12)(-0.5,-0.866)
\or\qbezier(-0.866,0.5)(-0.311,-0.179)(0,-1)
\or\qbezier(-0.866,0.5)(-0.146,-0.146)(0.5,-0.866)
\or\put(-0.866, 0.5){\line(5, -3){1.73}}
\fi
\or\ifcase#2\or\or\or\or\or\or
\or\qbezier(-1,0)(-0.516,-0.138)(-0.866,-0.5)
\or\qbezier(-1,0)(-0.45,-0.26)(-0.5,-0.866)
\or\qbezier(-1,0)(-0.327,-0.327)(0,-1)
\or\qbezier(-1,0)(-0.179,-0.311)(0.5,-0.866)
\or\qbezier(-1,0)(-0.0536,-0.2)(0.866,-0.5)
\fi
\or\ifcase#2\or\or\or\or\or\or\or
\or\qbezier(-0.866,-0.5)(-0.378,-0.378)(-0.5,-0.866)
\or\qbezier(-0.866,-0.5)(-0.26,-0.45)(0,-1)
\or\qbezier(-0.866,-0.5)(-0.12,-0.446)(0.5,-0.866)
\or\qbezier(-0.866,-0.5)(0,-0.359)(0.866,-0.5)
\fi
\or\ifcase#2\or\or\or\or\or\or\or\or
\or\qbezier(-0.5,-0.866)(-0.138,-0.516)(0,-1)
\or\qbezier(-0.5,-0.866)(0,-0.52)(0.5,-0.866)
\or\qbezier(-0.5,-0.866)(0.12,-0.446)(0.866,-0.5)
\fi
\or\ifcase#2\or\or\or\or\or\or\or\or\or
\or\qbezier(0,-1.)(0.138,-0.516)(0.5,-0.866)
\or\qbezier(0,-1.)(0.26,-0.45)(0.866,-0.5)
\fi
\or\ifcase#2\or\or\or\or\or\or\or\or\or\or
\or\qbezier(0.5,-0.866)(0.378,-0.378)(0.866,-0.5)
\fi\fi\fi\fi}
\def\FullChord[#1,#2]{
\Endpoint[#1]
\Endpoint[#2]
\Arc[#1]
\Arc[#2]
\Chord[#1,#2]
}
\def\EndChord[#1,#2]{
\Endpoint[#1]
\Endpoint[#2]
\Chord[#1,#2]
}
\def\Picture#1{
\begin{picture}(2,1)(-1,-0.167) 
#1
\end{picture}
}
\def\DottedChordDiagram[#1,#2]{
\Picture{\DottedCircle \FullChord[#1,#2]}
}
\def\ExtChord[#1,#2]{
\Endpoint[#1]\Endpoint[#2]
\thinlines
\ifnum#1>#2\ExtChord[#2,#1]
\else\ifnum#1<#2
\ifcase#1
\or\ifcase#2
\or\or\or\or\or\or\or\or\or
\qbezier[80](0,-1)(-0.1,-1.4)(0.25,-1.35)
\qbezier[80](0.25,-1.35)(1.35,-1.2)(1.35,0)
\qbezier[80](1.35,0)(1.35,0.95)(0.866,0.5)
\fi
\or\or\ifcase#2
\or\or\or\or\or\or\or\or\or\or\or
\qbezier[80](0,1)(-0.1,1.4)(0.25,1.35)
\qbezier[80](0.25,1.35)(1.35,1.2)(1.35,0)
\qbezier[80](1.35,0)(1.35,-0.95)(0.866,-0.5)
\fi
\or\ifcase#2
\or\or\or\or\or\or\or\or
\qbezier[80](-0.5,0.866)(-0.65,1.1)(-0.75,1.1)
\qbezier[60](-0.75,1.1)(-1.35,1.1)(-1.35,0)
\qbezier[60](-0.5,-0.866)(-0.65,-1.1)(-0.75,-1.1)
\qbezier[80](-0.75,-1.1)(-1.35,-1.1)(-1.35,0)
\fi\fi\fi\fi
}

\begin{document}
{\parindent0cm

\title{\bf Chromatic weight systems and the corresponding knot invariants}
\author{Jens Lieberum}
\maketitle

\begin{abstract}
{
\parindent0cm
We prove that chromatic weight systems, introduced by Chmutov, Duzhin and
Lando, can be expressed in terms of weight systems associated with direct
sums of the Lie algebras $\gl_n$ and $\so_n$.
As a consequence the Vassiliev invariants of knots
corresponding to the chromatic weight systems distinguish exactly the same
knots as a one variable specialisation $\V1$ of the
Homfly and Kauffman polynomial.
}
\end{abstract}

\section*{Introduction}

Vassiliev invariants of knots in $\R^3$ are in one-to-one correspondence with
weight systems, that is with linear forms on a vector space spanned by so-called
chord diagrams.
On one hand side, there exists a construction of weight systems using
reductive Lie algebras. On the other hand side, Chmutov, Duzhin and Lando introduced
chromatic weight systems without using Lie algebras.
This paper will give a positive and complete answer to the following natural question:

\begin{itemize}
\item Is there a relation between chromatic weight systems and weight systems 
associated to Lie algebras ?
\end{itemize}

We will also study the implications of this
relationship on the level of weight systems and on the level of knot invariants.

The paper is organized as follows.
From Section~\ref{algebraA} to Section~\ref{Wrficws} we prove that chromatic weight 
systems are linear combinations of weight systems associated with Lie algebras. 
For the proof we express a weight system called $\Wrxi$ in two different ways:

\begin{itemize}
\item The map $\Wrxi$ appears as the top coefficient of a polynomial expression
for weight systems coming from the Lie algebras $\gl_n^{\oplus r}$ 
(or $\so_n^{\oplus r}$)
and the $r$-th tensor power of the standard representation. This establishes
the connection to Lie algebras.

\item The map $\Wrxi$ can be expressed as a sum over colorings of the intersection graph
of a chord diagram using $r$ colors. This leads to the connection with chromatic weight
systems.
\end{itemize}

In Section~\ref{embed} we derive from the formula for $\Wrxi$ a formula for the value 
of the universal chromatic weight
system in terms of embeddings of certain trivalent diagrams into an oriented $2$-sphere.
In Section~\ref{other} we prove that chromatic weight systems are {\em not} linear
combinations of weight systems coming from direct sums of exceptional Lie algebras and
certain Lie superalgebras which completes our results on the level of weight 
systems.

From Section~\ref{algfknots} to Section~\ref{Homfly} we translate the results from
weight systems to Vassiliev invariants. 
We recall the relation between a weight system associated 
to $\gl_n$ and the Homfly polynomial $H$.
Then we extract from $H$ the part $\V1$ belonging to a certain chromatic weight system
$W$. The polynomial $\V1$ turns out to distinguish 
the same knots as all Vassiliev invariants
coming from chromatic weight systems. It is also a specialization of the Kauffman 
polynomial.

By \cite{Vo2} there exist weight systems that do not come 
from semisimple Lie superalgebras. This was an
inspiration for our research, as well as the need of a good combinatorial 
understanding of the known weight systems.

\section*{Acknowledgement}

The starting-point for this work was in Bonn where I was guest at the 
'Max-Planck-Institut f\"{u}r Mathematik' and where C.\ F.\ B\"{o}digheimer told me about
chromatic weight systems. Most of the work for this paper was carried out at the
Universit\'{e} Louis Pasteur at Strasbourg. 
I would like to express my deep gratitude to C.\ Kassel for his help, especially for 
giving me important literature references and for proposing many small and great changes
in earlier versions of this paper.
I thank the German Academic Exchange Service for financial support 
(Doktorandenstipendium HSP II/AUFE).

\section{The bialgebra of chord diagrams}\label{algebraA}

Consider finite unoriented graphs $D$ (multiple edges are allowed) with the following
properties:

1) $D$ is equipped with an embedding of an oriented circle $S^1\subset D$.

2) Every vertex has valency three.

3) The vertices outside of the image of $S^1$ are equipped with a cyclic order 
on the three arriving edges.

Graphs with these properties are said to be equivalent if there exists a
ho\-meo\-mor\-phism of graphs that respects the embeddings of
$S^1$ with orientation and preserves the cyclic order of the vertices outside of $S^1$. 
An equivalence class of these graphs will be called a {\em trivalent diagram}.
A trivalent diagram $D$ will be called a {\em chord diagram} if all vertices lie on
$S^1\subset D$. We define the degree of $D$ by

\begin{equation}
\deg D:=\frac{1}{2}\#\{\mbox{trivalent vertices of $D$}\}\in \N.
\end{equation}

In our pictures of trivalent diagrams the embedded oriented circle $S^1$ 
is drawn with a thick
line and is always oriented counterclockwise. The rest of the diagram is drawn
in the interior of the disk bounded by $S^1$. The cyclic order of the trivalent 
vertices outside of $S^1$ is counterclockwise. Vertices with valency four in our 
pictures do not correspond to vertices of the graph. They only appear
because we draw the diagrams in the plane. We use dots to indicate that we only show
a part of a trivalent diagram. 

Let $R$ be a field of characteristic 0. The symbol $\otimes$ will always 
denote the tensor product over $R$.
If $B$ is a graded $R$-\vs, we
will denote the space of homogeneous elements of degree $d$ by $B_d$ and 
the space spanned by all homogeneous elements of degree $\leq d$ ($\geq d$)
by $B_{\leq d}$ ($B_\geq d$).

\begin{defi}\label{algA}
Let $A$ be the graded $R$-vector space with trivalent graphs as generators
and with the so-called (STU)-relation from Picture~(\ref{sturel}) for all triples of 
diagrams that only differ like the three diagrams in this picture.
\end{defi}

\begin{center}
\begin{equation}\label{sturel}
\Picture{
\DottedCircle
\Arc[8]
\Arc[9]
\Arc[10]
\Endpoint[9]
\thinlines
\put(0.0,-1.0){\line(0,1){0.6}}
\put(0.0,-0.4){\line(1,1){0.5}}
\put(0.0,-0.4){\line(-1,1){0.5}}
}\quad = \quad
\Picture{
\DottedCircle
\Arc[8]
\Arc[9]
\Arc[10]
\Endpoint[8]
\Endpoint[10]
\thinlines
\put(0.5,-0.866){\line(0,1){0.966}}
\put(-0.5,-0.866){\line(0,1){0.966}}
} \quad - \quad 
\Picture{
\DottedCircle
\Arc[8]
\Arc[9]
\Arc[10]
\Endpoint[8]
\Endpoint[10]
\thinlines
\put(0.5,-0.866){\line(-1,1){0.966}}
\put(-0.5,-0.866){\line(1,1){0.966}}
}
\end{equation}
\\[12pt]
\end{center}

The implicit statement in Definition~\ref{algA} that $A$ is a graded \vs{} holds
because the (STU)-relation is homogeneous with respect to the degree of trivalent
diagrams.
The (STU)-relation allows one to write every element of $A$ as a linear combination of
chord diagrams. Applying the (STU)-relation two times gives the so called (4T)-relation
($=$ four-term relation) from Picture~(\ref{ftrel}) for chord diagrams.

\begin{center}
\begin{equation}\label{ftrel}
\Picture{\DottedCircle\FullChord[1,5]\FullChord[9,0]\Arc[6]}
\quad - \quad
\Picture{\DottedCircle\FullChord[0,5]\FullChord[9,1]\Arc[6]}
\quad = \quad
\Picture{\DottedCircle\FullChord[1,5]\FullChord[9,6]\Arc[0]}
\quad - \quad
\Picture{\DottedCircle\FullChord[1,6]\FullChord[9,5]\Arc[0]}
\end{equation}\\[12pt]
\end{center}

As a consequence there exists a surjective morphism from the \vs{} generated
by chord diagrams and (4T)-relations to $A$ induced by the inclusion of sets of diagrams.
This map is an isomorphism (see \cite{BN1}, Theorem 6). 

Now we describe the Hopf algebra structure of $A$.
Given two trivalent diagrams $D_1$ and $D_2$, cut their oriented circles somewhere
outside of a vertex and glue the resulting oriented intervals to form one new oriented
circle.
This operation will be called the connected sum of $D_1$ and $D_2$. It can be proved that
the connected sum is well defined for classes of diagrams in $A$, so it turns $A$ into a
commutative graded algebra with unit $S^1$.

For a trivalent diagram $D$ define

\begin{equation}\label{deltar}
\Delta_r(D):=\sum_{D=D_1 \cup\ldots\cup D_r } D_1 \otimes\ldots \otimes D_r
\end{equation}

where in the case of $\deg D>0$ the sum is taken over all functions
$\varphi:\pi_0(D\setminus S^1)\longrightarrow \{1,\ldots, r\}$ 
and the diagrams $D_1,\ldots,D_r$ in the sum are defined by

$$
D_i=D\setminusÊ\bigcup_{E:\varphi(E)\not=i} E.
$$

For $D=S^1$ the notation in Formula~(\ref{deltar}) will mean that
$\Delta_r(S^1)=S^1\otimes\ldots\otimes S^1$.

The map $\Delta:=\Delta_2$ turns $A$ into a commutative, cocommutative bialgebra and
we have $\Delta_r=(\Delta_{r-1}\otimes id)\circ \Delta$. 
The counit $\epsilon$ is determined by 

\begin{equation}\label{counit}
\epsilon(D)=
\left\{
\begin{array}{ll}
1 & \mbox{if $D=S^1$},\\
0 & \mbox{if $\deg D>0$.}
\end{array}
\right.
\end{equation}

We also have $\dim A_0=1$ and $\dim A_i<\infty$ for all $i$. 
A structure theorem (\cite{MM65}) implies that 
the graded bialgebra $A$ is a Hopf algebra isomorphic to the algebra of polynomials on 
the graded \vs{}

$$
P(A):=\{a\in A\ \vert \ \Delta(a)=a\otimes 1+1\otimes a\}
$$

of primitive elements of $A$.
By the definition of $\Delta$ it is clear that

\begin{equation}\label{primel}
{\cal M}:=\left\{
{\mbox{trivalent diagrams $D$ such that}\atop 
\mbox{$D\setminus S^1$ is connected}}
\right\}\subset P(A).
\end{equation}

On the other hand, the (STU)-relation allows one to write every trivalent diagram as a
linear combination of products of elements from ${\cal M}$. This implies that 
$\span {\cal M}=P(A)$.
We do not need the so-called Chinese characters (see \cite{BN1}) for 
this purpose. Useful applications of these diagrams different to the ones in 
\cite{BN1} are a decomposition of 
$A_d$ into a direct sum of eigenspaces of the
cabling operations (\cite{KSA}) and a proof that slight modifications of the 
definition of $P(A)$ are possible (Proposition~1.1 of \cite{Vo2}). 
We will use the last application implicitely in Section~\ref{other}.

\section{Weight systems}

An open question concerning the structure of $A$ is the determination of (the 
asymptotic behaviour of) the sequence $(\dim P(A)_i)_{i\geq 10}$. The best we can do 
in higher degrees in order to
obtain some information from a linear combination of classes of trivalent diagrams 
in $A$ is to apply an element of the dual space of $A$ to this linear combination.

\begin{defi}
Call a linear map $A_d\longrightarrow R$ a weight system of degree $d$, a linear map 
$A\longrightarrow R$ a weight system, and an algebra homomorphism $A\longrightarrow R$ a 
multiplicative weight system.
\end{defi}

The dual space $A^*\cong \prod_{i=0}^\infty A_i^*$ of all weight systems is an algebra
and the graded dual space $\bigoplus_{i=0}^\infty A_i^*$  generated by all weight 
systems of finite degree is a graded Hopf algebra. The algebra $A^*$ becomes a Hopf 
algebra in a
category with a completed tensor product and the multiplicative weight systems are
the group-like elements of this Hopf algebra.
Instead of $R$ we may take some $R$-\vs{} $M$ as the space of values of weight 
systems. Composition with elements of $M^*$ gives us weight systems with values in $R$.
A weight system $w$ with values in a finite dimensional \vs{} $M$ is said to be a 
linear combination of other weight systems with values in $R$ if this is true for 
$p\circ w$ for all $p\in M^*$.

Weight systems can be obtained from Lie superalgebras. 
A Lie superalgebra $\L=\L_{\overline{0}}\oplus\L_{\overline{1}}$ is a 
$\Z/(2)$-graded vector space with a graded bilinear map $\L\times\L\longrightarrow\L$ 
satisfying an antisymmetry
property and a Jacoby identity. 
Lie algebras are Lie superalgebras with $\L_{\overline{1}}=\{0\}$. We will consider
the case $\L_{\overline{1}}\not=\{0\}$ only in this section and in Section~\ref{other}. 
Let $(\ ,\ )$ be a consistent, supersymmetric 
$\ad$-invariant, non-degenerate bilinear form on the Lie superalgebra $\L$ 
(see \cite{Kac} for the terminology). 
An example of such a bilinear form is the 
Killing form for a semisimple Lie algebra.
We choose two homogeneous bases $e_i$ and $e_j'$ of $\L$ such that
$(e_j',e_i)=\delta_{ij}$ ($i,j=1,\ldots,\dim \L$) and define 
$\omega:=\sum_i e_i\otimes e_i'$. The definition of 
$\omega$ does not depend on the choice of the bases and $\omega$ is supersymmetric,
$\ad$-invariant and of degree $0$. 
We call $\omega$ a Casimir element.
Let $\rho:\L\longrightarrow \End M$ be a representation of 
$\L$ in a $\Z/(2)$-graded, finite-dimensional \vs{} $M$. 
To $(\rho\otimes\rho)(\omega)$ we will associate a weight system as follows:

Let $D$ be a chord diagram of degree $n$. Then $D$ looks as in 
Picture~(\ref{permutcd})
for a non-unique permutation $\pi_D\in S_{2n}$ where in the box named $\pi_D$
the top chord-end number $i$ is connected to the bottom chord-end number $\pi_D(i)$.

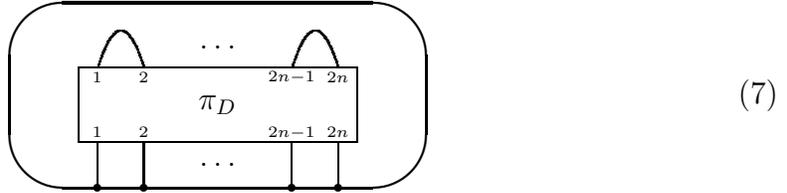
\begin{center}
{\unitlength=35pt
\begin{equation}\label{permutcd}
\begin{picture}(5, 1)(-2.5, -0.08)
\thicklines
\put(0,0){\oval(4.5,2)}
\thinlines
\put(-1.5,-0.5){\framebox(3.0,0.8){$\pi_D$}}
\put(-1.3,-1){\line(0,1){0.5}}
\put(-0.8,-1){\line(0,1){0.5}}
\put(0.8,-1){\line(0,1){0.5}}
\put(1.3,-1){\line(0,1){0.5}}
\put(-1,-1){\makebox(2,0.5){$\ldots$}}
\qbezier(-1.3,0.3)(-1.05,1.1)(-0.8,0.3)
\qbezier(1.3,0.3)(1.05,1.1)(0.8,0.3)
\put(-1,0.3){\makebox(2,0.45){$\ldots$}}
\put(-1.5,0.1){\makebox(0.4,0.2){$\scriptscriptstyle 1$}}
\put(-1.0,0.1){\makebox(0.4,0.2){$\scriptscriptstyle 2$}}
\put(0.6,0.1){\makebox(0.4,0.2){$\scriptscriptstyle 2n-1$}}
\put(1.1,0.1){\makebox(0.4,0.2){$\scriptscriptstyle 2n$}}
\put(-1.5,-0.5){\makebox(0.4,0.2){$\scriptscriptstyle 1$}}
\put(-1.0,-0.5){\makebox(0.4,0.2){$\scriptscriptstyle 2$}}
\put(0.6,-0.5){\makebox(0.4,0.2){$\scriptscriptstyle 2n-1$}}
\put(1.1,-0.5){\makebox(0.4,0.2){$\scriptscriptstyle 2n$}}
\put(1.3,-1){\circle*{0.08}}
\put(0.8,-1){\circle*{0.08}}
\put(-0.8,-1){\circle*{0.08}}
\put(-1.3,-1){\circle*{0.08}}
\end{picture}
\end{equation}\nopagebreak\vspace*{21pt}
}
\end{center}

The symmetric group $S_{2n}$ acts on $(\End M)^{\otimes 2n}$ from the left by 
superpermutation of the tensor factors:

\begin{equation}
\pi\cdot \varphi_1\otimes\ldots\otimes\varphi_{2n}:=
(-1)^{\sigma(\pi,\deg \varphi_1,\ldots,\deg \varphi_{2n})}
\varphi_{\pi^{-1}(1)}\otimes\ldots\otimes\varphi_{\pi^{-1}(2n)}
\end{equation}

where

$$
\sigma(\pi,a_1,\ldots,a_{2n}):=\#\left\{
(i,j)\mid 1\leq i<j\leq 2n,\pi(i)>\pi(j),a_i=a_j=\overline{1}\right\}.
$$

Let $\gamma_{2n}:\End M^{\otimes 2n}\longrightarrow \End M$ be the composition of $2n$ 
endomorphisms and define

\begin{equation}\label{Lieweight}
W_{\L,\omega,\rho}(D):=\left\{
\begin{array}{cl}
\str\left(\gamma_{2n}\left(\left(\pi_D\cdot\left(\rho^{\otimes 2n}(\omega^{\otimes n})
\right)\right)\right)\right) & \mbox{if $\deg D>0$,}\\
\sdim M:=\dim M_{\overline{0}}-\dim M_{\overline{1}} & \mbox{if $D=S^1$}
\end{array}
\right.
\end{equation}

where $\str$ denotes the supertrace.

\begin{lemma}\label{welldef}
Formula (\ref{Lieweight}) defines a weight system 
$W_{\L,\omega,\rho}:A\longrightarrow R$.
\end{lemma}

For more useful details, extensions or generalisations of this construction of weight 
systems and for a proof of Lemma~\ref{welldef} 
see Section~6 of \cite{Vo2}, \cite{Va1} and Section~2.4 of \cite{BN1}.

Observe that for $a\in R$ we have

\begin{equation}\label{Casinorm}
W_{\L,a\omega,\rho}(D)=a^{\deg D} W_{\L,\omega,\rho}(D).
\end{equation}

Define a contraction map $\kappa_{2n}:(M\otimes M^*)^{\otimes 2n}
\longrightarrow R$ by

\begin{equation}
\kappa_l(m_1 \otimes \varphi_2\otimes m_2 \otimes \ldots 
\otimes \varphi_{2n}\otimes m_{2n}\otimes \varphi_1):=(-1)^{\deg m_1}
\prod_{i=1}^{2n} \varphi_i(m_i).
\end{equation}

Let $\iota_M:\End M\longrightarrow M\otimes M^*$ be the canonical isomorphism. Then 
we can restate Formula~(\ref{Lieweight}) for $\deg D>0$ as

\begin{equation}\label{Fbigsum}
W_{\L,\omega,\rho}(D)=\left(\kappa_{2n}\circ \iota_M^{\otimes 2n}\right)\left(\pi_D\cdot\left( 
\rho^{\otimes 2n}(\omega^{\otimes n})\right)\right).
\end{equation}

\section{Chromatic weight systems}\label{Cws}

Now we review some facts from \cite{CDL2} and \cite{CDL3} concerning chromatic weight
systems.
Consider a chord diagram $D$. We call $\pi_0(D\setminus S^1)$ the set of chords 
of $D$. We say that the chords $a$ and $b$ intersect if we meet the endpoints of $a$
and $b$ in alternating order when we travel around the circle of~$D$.

\begin{defi}
The intersection graph $\Gamma(D)$ of a chord diagram $D$ is a graph 
whose vertices are in bijection with the chords of $D$ and two
vertices are connected by an edge exactly if the corresponding chords of $D$ intersect.
\end{defi}
 
Now consider graphs $G$ without orientation, without loops and without multiple edges
and let $\MO (G)$ denote the set of vertices of $G$.
A weighted graph $G$ is a graph together with a function $w : \MO (G)\longrightarrow \N$. 
Define the degree of $G$ as the sum of the weights of all vertices.
Given an edge $e=\edge[v_1,v_2]$ of a weighted graph $G$, we can define two new
weighted graphs $G_e'$ and $G_e''$:

\begin{itemize}
\item $G_e'$ is obtained from $G$ by removing the edge $e$.

\item $G_e''$ is obtained from $G_e'$ by first contracting the vertices $v_1$ and $v_2$
to one new vertex with the weight $w(v_1)+w(v_2)$ and then
by replacing multiple edges by unique edges.
\end{itemize}

\begin{defi}
Let ${\cal C}$ be the graded $R$-vector space with weighted graphs including
the empty graph as generators 
and with the 'chromatic' relations $G=G_e' - G_e''$  
for all edges $e$.
\footnote{The relations $G=G_e' + G_e''$ of \cite{CDL3} lead (up to signs) to
the same, but would cause slightly more complicated formulas here.
The chromatic relation should be compared to the (STU)-relation.}
\end{defi}

We have a map from chord diagrams to weighted graphs such that a chord diagram $D$ is
mapped to $\Gamma(D)$ with the weight of all vertices being one. This map
induces a linear map 

\begin{equation}
\OverGam:A\longrightarrow {\cal C}
\end{equation}

of graded modules and there exists a structure of a
Hopf algebra on ${\cal C}$ that turns this map into a morphism of Hopf algebras.
The product of two weighted graphs is their disjoint
union and the coproduct is described in \cite{CDL3}.
A consequence of the chromatic relations is that ${\cal C}$ is spanned by products of
single vertices where the weights may be chosen freely. This means that there exists a 
surjective morphism of graded algebras 

\begin{equation}
R[s_1,s_2,\ldots]\longrightarrow {\cal C}\ \ , \ \ s_d\mapsto
\unitlength=1ex
\begin{picture}(2.0,2.5)(0.0,0.5)
\put(0.7,-0.05){\circle*{0.5}}
\put(0.3,0.95){\makebox(1,1){$d$}}
\end{picture}
\ \ , \ \ (\deg s_d:=d)
\end{equation}

sending $s_d$ to the graph consisting of a single vertex with weight $d$.
By Theorem~2 of \cite{CDL3} this map is an isomorphism. 

\begin{defi}
A weight system that factors through ${\cal C}$ will be called chromatic 
and $\OverGam$ will be called the universal chromatic weight system.
\end{defi}

The space of chromatic weight systems is a subalgebra of $A^*$ isomorphic to 
${\cal C}^*$ because $\OverGam$ is
a surjective morphism of coalgebras. 

We have described two types of weight systems that are
not obviously related to each other.
One of our main goals is to prove that chromatic weight systems come from 
Lie algebras.

\begin{theorem}\label{linearcomb}
Every chromatic weight system of degree $d$ is a linear combination of weight systems 
associated with
the Lie algebras $\gl_n^{\oplus d}$ and also a linear combination of weight systems 
associated with the Lie algebras $\so_n^{\oplus d}$.
\end{theorem}

For the previous theorem we can use for both Lie algebras the $d$-th tensor power 
of the 
defining representation. The Casimir element and $n$ vary for the weight systems used
in the linear combination.
We prove Theorem~\ref{linearcomb} in Section~\ref{Wrficws}.

\section{Weight systems coming from $\gl_n$}\label{glnws}

For a chord diagram $D$ define

$$
W(D):=\left\{
\begin{array}{ll}
0 & \mbox{if $D$ has two intersecting chords,}\\
1 & \mbox{otherwise.}
\end{array}
\right.
$$

It will follow from Lemma~\ref{Wgltau} and is easy to see directly 
that $W$ induces a multiplicative weight system $W:A\longrightarrow R$. 
In each finite degree we want to obtain $W$ as a linear combination of weight systems
coming from the Lie algebras $\gl_n$ and $\so_n$.
We start with a purely combinatorial lemma.
It is formulated so as to be useful for the proofs of Lemma~\ref{Wgltau},
Lemma~\ref{Wso}, and also for a possible version of Lemma~\ref{Wso} for the
Lie algebra $\spli_n$ (see Exercise~6.33 and Exercise~6.34 of \cite{BN1} or write down
an explicit formula for a Casimir element of $\spli_n$).

\begin{lemma}\label{combi}
Let $D$ be a chord diagram with $d$ chords. Forget the orientation of $S^1$. If we
replace each chord by one of the three possibilities from Picture~(\ref{threepos}),
then we get $k\leq d+1$ circles. We get $d+1$ circles exactly 
if no chords of $D$ intersect and every chord is replaced by Possibility (1) of
Picture~(\ref{threepos}).

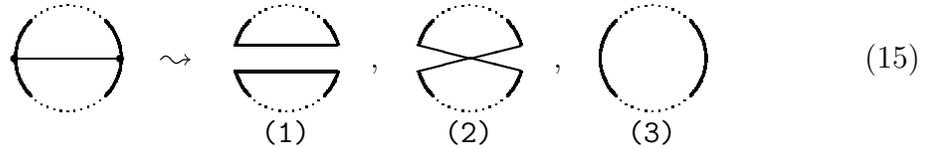
\begin{center}
\begin{equation}\label{threepos}
\Picture{\DottedCircle\FullChord[6,0]\Arc[1]\Arc[5]
\Arc[7]\Arc[11]}
\quad \leadsto \quad
\Picture{
\Arc[1]\Arc[5]
\Arc[7]\Arc[11]
\thinlines
\DottedArc[2]
\DottedArc[3]
\DottedArc[4]
\DottedArc[8]
\DottedArc[9]
\DottedArc[10]
\thicklines
\put(-0.975,-0.25){\line(1,0){1.95}}
\put(-0.975,0.25){\line(1,0){1.95}}
\put(-1,-1.7){\makebox(2,0.6){\tt (1)}}
}
\quad , \quad
\Picture{
\Arc[1]\Arc[5]
\Arc[7]\Arc[11]
\thinlines
\DottedArc[2]
\DottedArc[3]
\DottedArc[4]
\DottedArc[8]
\DottedArc[9]
\DottedArc[10]
\thicklines
\put(-0.97,-0.25){\line(4,1){1.95}}
\put(-0.97,0.25){\line(4,-1){1.95}}
\put(-1,-1.7){\makebox(2,0.6){\tt (2)}}
}
\quad , \quad
\Picture{\DottedCircle\Arc[0]\Arc[1]\Arc[6]
\Arc[5]\Arc[7]\Arc[11]
\put(-1,-1.7){\makebox(2,0.6){\tt (3)}}
}
\end{equation}\nopagebreak\vspace*{24pt}
\end{center}
\end{lemma}

{\bf Proof:}
The lemma can be proved by induction on the number of chords
using chord diagrams on some number of circles.
$\Box$

\medskip

We recall the graphical calculus for a weight system associated to $\gl_n$ from 
\cite{BN1}.
Let $M$ be an $n$-dimensional vector space with basis $e_i$ ($i=1,\ldots,n$).
Let $e_{ij}$ denote the standard basis of $\End M$ and let 
$\tau:=\id:\gl_n\longrightarrow \End M$ be
the defining representation of the Lie algebra $\gl_n$. 
Define the $\ad$-invariant, symmetric, 
non-degenerate bilinear form $(\ ,\ )$ by $(a,b):=\tr(\tau(a)\tau(b))$. The 
corresponding Casimir element is

\begin{equation}\label{casimirX}
X=\sum_{i,j=1}^n e_{ij}\otimes e_{ji}
\in \gl_n\otimes \gl_n.
\end{equation}

As a linear map $X$ is the permutation of the tensor factors of $M\otimes M$.
Now recall the expression $W_{\gl_n,X,\tau}$ from Formula~(\ref{Fbigsum}). 
Let $e^i$ be the basis of $M^*$ dual to $e_i$ ($i=1,\ldots,n$).
Then 

\begin{equation}
(\iota_M\otimes\iota_M)\circ (\tau\otimes\tau)(X)
=\sum_{i,j=1}^n e_i\otimes e^j\otimes e_j\otimes e^i.
\end{equation}

This tensor appears in a 'big sum expansion' of $W_{\gl_n,X,\tau}(D)$ with 
$n^{4\deg D}$ terms as indicated by Picture~(\ref{bigsum}).

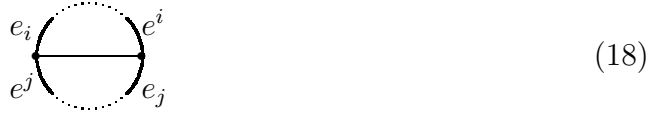
\begin{center}
\begin{equation}\label{bigsum}
\begin{picture}(2,1)(-1,-0.167)
\DottedCircle\FullChord[6,0]
\Arc[1]\Arc[5]\Arc[7]\Arc[11]
\put(-1.5,0.4){\makebox{$e_i$}}
\put(1.0,0.4){\makebox{$e^i$}}
\put(-1.5,-0.8){\makebox{$e^j$}}
\put(1.0,-0.8){\makebox{$e_j$}}
\end{picture}
\end{equation}\nopagebreak\vspace*{12pt}
\end{center}

A term with two vectors $e_{i_\nu}$ and $e^{i_\mu}$ that label an interval on the 
oriented circle between two neighbored chord-ends
can only give a contribution to the 'big sum' if $i_\nu=i_\mu$. 
This and Picture~(\ref{bigsum}) implies that

\begin{equation}\label{WglnXtau}
W_{\gl_n, X, \tau}(D)=n^{c(D)}
\end{equation}

where $c(D)$ denotes the number of circles that appear when all chords of $D$ are
replaced by Possibility~(1) of Picture~(\ref{threepos}). We summarize:

\begin{lemma}\label{Wgltau}
The formula 

$$
W_{\gl}(D)=t^{c(D)-1}
$$

for chord diagrams $D$
defines a homomorphism of algebras $W_{\gl}:A\longrightarrow R[t]$
with the property

$$
W_{\gl_n,X,\tau}(D)=n W_{\gl}(D)(n)
$$

for all $n\geq 1$.
If $D$ has degree $d$, then $\deg(W_{\gl}(D))\leq d$ and the coefficient
of $t^d$ is $W(D)$.\footnote{The constant term in $W_{\gl}$ is the 
weight system of the Alexander-Conway polynomial and $W_{\gl}$ can also be seen to
be well-defined for geometric reasons (see \cite{BNG}, Section~3).}
\end{lemma}
{\bf Proof:} The map $W_{\gl}$ is well defined because in Formula 
(\ref{WglnXtau}) infinitely many choices of $n$ are possible. It is 
multiplicative with respect to the connected sum of chord diagrams because $c$ is
additive with respect to the disjoint union of chord diagrams.
By Lemma~\ref{combi} we have $c(D)\leq d+1$ for all chord diagrams $D$ and 
$c(D)=d+1$ exactly if no chords of $D$ intersect.
This completes the proof.$\Box$

\medskip

We can formulate a version of Lemma~\ref{Wgltau} also for the Lie algebras
$\so_n$. 

\begin{lemma}\label{Wso}
Let $\tau$ be the defining representation of $\so_n$.
We use a bi\-li\-near form $(\ ,\ )$ defined by $(a,b):=\frac{1}{2}\tr(\tau(a)\tau(b))$
and denote the associated Casimir element by $\omega$.
Then there exists a morphism of $R$-algebras $W_{\so}:A\longrightarrow
R[t]$ such that for all $n\geq 2$ we have 

$$
W_{\so_n,\omega,\tau}(D)=n W_{\so}(D)(n).
$$

If $D$ has degree $d$, then $\deg(W_{\so}(D))\leq d$ and the coefficient
of $t^d$ is $W(D)$.
\end{lemma}

A proof of Lemma~\ref{Wso} follows from Formula~(29) of \cite{BN1} and
Lemma~\ref{combi}.

\section{Products of weight systems}\label{products}

For a chord diagram $D$ define

\begin{equation}\label{defwrxi}
\Wrxi (D):=
\sum_{D=D_1\cup\ldots\cup D_r}
x_1^{\deg D_1} \ldots x_r^{\deg D_r}W(D_1\ldots D_r)\in R[x_1,\ldots,x_r].
\end{equation}

It will follow from Proposition~\ref{WLVr} 
and is easy to see directly that the map $\Wrxi$ is well defined. In each finite 
degree we want to obtain $\Wrxi$ as a linear combination of weight systems coming from
the Lie algebras $\gl_n^{\oplus r}$ and $\so_n^{\oplus r}$. 
In the sequel we will suppress some indices and
hope that this does not lead to confusion.

Let $\L_n$ be one of the Lie algebras $\gl_n$ or $\so_n$, let $\tau$ be the defining 
representation and $\omega$ be the Casimir element from Lemma \ref{Wgltau} or
Lemma \ref{Wso}.
Let $a_1,\ldots,a_r\in R$, define
$\Omega:=(a_1 \omega,\ldots,a_r \omega)
\in (\L_n\otimes \L_n)^{\oplus r}\subset (\Lnr)\otimes (\Lnr)$
and let the $\Lnr$-representation 
$\tau^{\otimes r}$ be the $r$-fold 
tensor power of the standard representation~$\tau$.

\begin{prop}\label{WLVr}
For $\L=\gl$ or $\so$ and $r\geq 1$ there exists a morphism of $R$-algebras 

$$
W_{\L,r}:A\longrightarrow \sympol[t]
$$

such that for all choices of $a_1,\ldots,a_r\in R, n\geq 2$ we have

$$
W_{\Lnr, \Omega, \tau^{\otimes r}}(D)=n^r W_{\L, r}(D)(a_1,\ldots,a_r,n).
$$

If $D$ has degree $d$, then the degree in the indeterminate $t$ of $W_{\L, r}(D)$
is $\leq d$ and the coefficient of $t^d$ is $\Wrxi(D)$.
\footnote{The constant term of $W_{\gl,r}$ can be expressed in terms of immanent
weight systems (see Section~6 of \cite{BNG} and Section~5 of \cite{KSA}).}
\end{prop}
{\bf Proof:}
Formula~(\ref{Casinorm}) and Exercise~6.33 of \cite{BN1} give us the formula

$$
W_{\Lnr, \Omega, \tau^{\otimes r}}(D)=\sum_{D=D_1\cup\ldots\cup D_r}
a_1^{\deg D_1} W_{\L_n,\omega,\tau}(D_1) \ldots a_r^{\deg D_r}
W_{\L_n,\omega,\tau}(D_r).
$$

This formula is valid for all choices of $a_1,\ldots,a_r\in R$ and $n\geq 2$ and
$W_{\L_n,\omega,\tau}(D_i)$ is a polynomial in $n$ with vanishing constant term
(see Lemma~\ref{Wgltau} and Lemma~\ref{Wso}), so there
exists a unique weight system $W_{\L,r}$ with the defining property from the 
proposition.
We again use the Lemmas~\ref{Wgltau} and \ref{Wso} to see that $\Wrxi(D)$ is the
coefficient of 
$t^d$ in $W_{\L,r}(D)$ ($\deg D=d$) and to see that $W_{\L, r}$ is a morphism of
algebras.
We have $W_{\L,r}(a)\in\sympol[t]$ for all $a\in A$ because $A$ is 
cocommutative.
$\Box$

\medskip

Recall that the algebra $R[x_1,\ldots,x_r]^{S_r}$ is a polynomial algebra in the 
polynomials

\begin{equation}
G_{i,r}:=x_1^i+\ldots+x_r^i\ \ (i=1,\ldots,r).
\end{equation}

In order to be as explicit as possible in the proof of Theorem~\ref{linearcomb} we
formulate the next lemma. 

\begin{lemma}\label{coefffkt}
The coefficient functions of monomials in $G_{1,r},\ldots,G_{r,r},t$ 
from $W_{\L, r}:A_d\longrightarrow R[G_{1,r},\ldots,G_{r,r}]_d [t]$ 
are linear combinations of weight systems associated to the Lie algebras $\Lnr$.
\end{lemma}

{\bf Proof:}
We give a short constructive proof and we do not care about the number of
terms in the linear combination:

Let $P_i$ ($i=2,\ldots,r+d+2$) be the unique polynomial of degree $r+d$, such that 
$P_i(j)=\delta_{ij}$ ($j=2,\ldots,r+d+2$). In order to expand the result in
the desired basis we choose some $R$-linear projection 
$\pi:R[x_1,\ldots,x_r]\longrightarrow R[G_{1,r},\ldots,G_{r,r}]_d$.
Then for chord diagrams $D$ of degree $d$ we have

$$
W_{\L, r}(D)=\sum_{a_1,\ldots,a_r,n=2}^{r+d+2}
W_{\Lnr,(a_1 \omega,\ldots, a_r \omega),\tau^{\otimes r}}(D)
\pi(P_{a_1}(x_1)\ldots P_{a_r}(x_r))
P_n(t)/n^r
$$

and the coefficients of monomials in $G_{1,r},\ldots,G_{r,r},t$ can be read 
from the right hand side as linear combinations of weight systems coming
from the Lie algebras~$\Lnr$.
$\Box$

\section{Chromatic weight systems and $\Wrxi$}\label{Wrficws}

Let $f_i : {\N}\longrightarrow R$ ($i=1,\ldots, r$) be functions.
Define $\eta_i:A\longrightarrow A$ by $\eta_i(a):=f_i(\deg a) a$.
We will be interested in the weight system

\begin{equation}\label{defwrfi}
\Wrfi :=W^{\otimes r} \circ (\eta_1\otimes\ldots\otimes \eta_r)
\circ \Delta_r.
\end{equation}

Let $a_i\in R$ ($i=1,\ldots,r$). For the special
choice $f_i(d)=a_i^d$ we have 
$\Wrfi(D)=\Wrxi(D)(a_1,\ldots,a_r)$ because of Formula~(\ref{defwrxi}). 
It will soon turn out that other choices of the maps 
$f_i$ will not lead to greater generality.
Let us now introduce some terminology to describe the combinatorics underlying 
Formula~(\ref{defwrfi}).

\begin{defi}
Let $G$ be a graph with vertex set $\MO(G)$. We call a function $c:\MO(G)\longrightarrow
\{1,\ldots, r\}$ a vertex coloring of $G$ with $r$ colors if for each edge 
$\edge[v_1,v_2]$ we have $c(v_1)\not=c(v_2)$. Define

$$
C_r(G):=\left\{
\begin{array}{cl}
\left\{
c:\MO(G)\longrightarrow\{1,\ldots, r\}\ \vert \ \mbox{$c$ is a vertex coloring}
\right\} & \mbox{if $G\not=\emptyset$},\\
\{\id:\{0\}\longrightarrow\{0\}\} & \mbox{if $G=\emptyset$}.
\end{array}
\right.
$$
\end{defi}

Notice that $W(D)=\# C_1(\Gamma(D))$.
The next proposition will be responsible for the connection between $\Wrfi$ and
chromatic weight systems.

\begin{prop}\label{colorsum1} The following formula holds for all chord diagrams $D$:
$$
\Wrfi (D)=\sum_{c\in C_r(\Gamma(D))}
\prod_{i=1}^r f_i(\# c^{-1}(i)).
$$
\end{prop}
{\bf Proof:}
By definition we have for a chord diagram $D$:
$$
\Wrfi (D)=\sum_{c:\MO(\Gamma(D))\longrightarrow \{1,\ldots,r\}}
\prod_{i=1}^r f_i(\# c^{-1}(i)) W(D_{c,i}).
$$

where $D_{c,i}$ denotes the subdiagram of $D$ with $\{a\ \vert\ c(a)=i\}$ as the set
of chords.
But

$$
\prod_{i=1}^r W(D_{c,i})=\left\{
\begin{array}{ll}
1 & \mbox{if $c$ is a vertex coloring of $\Gamma(D)$,}\\
0 & \mbox{otherwise.}
\end{array}
\right.
$$
$\Box$

\medskip

Let $G$ be a weighted graph with weight function $w$ and define

\begin{equation}\label{defwrfit}
\Wrfit(G):=\sum_{c\in C_r(G)}
\prod_{i=1}^r f_i \left(\sum_{v\in \MO(G)} \delta_{i,c(v)}w(v)\right).
\end{equation}

\begin{lemma}\label{l:colorsum2}
Formula~(\ref{defwrfit}) defines a linear map $\Wrfit:{\cal C}\longrightarrow R$. 
\end{lemma}
{\bf Proof:}
Let $e=\edge[v_1,v_2]$ be an edge of $G$ and recall the definition of $G_e'$ and 
$G_e''$ from Section~\ref{Cws}.
We can identify the sets 

$$
\left\{c \in C_r(G_e')\ \vert\ c(v_1)\not=c(v_2)\right\}\mbox{ and }
\left\{c \in C_r(G_e')\ \vert\ c(v_1)=c(v_2)\right\}
$$

with $C_r(G)$ and $C_r(G_e'')$ respectively. 
Let $c \in C_r(G_e')$ with $c(v_1)=c(v_2)$ and
denote by $\overline{c}$ the corresponding element of $C_r(G_e'')$.
Denote the weight functions of $G_e'$ and $G_e''$ by $w$ and $\overline{w}$ respectively.
Then we have for $i=1,\ldots, r$:

$$
\sum_{v\in \MO(G_e')} \delta_{i,c(v)}w(v)=\sum_{v\in \MO(G_e'')} 
\delta_{i,\overline{c}(v)}\overline{w}(v),
$$

because in the definition of $G_e''$ the new contracted vertex has the weight 
$w(v_1)+w(v_2)$. A similar statement holds for $c \in C_r(G_e')$ with 
$c(v_1)\not=c(v_2)$ and the corresponding element $\overline{c}\in C_r(G)$.
This implies

$$
\Wrfit(G)=\Wrfit(G_e')-\Wrfit(G_e'')
$$

and completes the proof.$\Box$

\medskip

We have $\Wrfit(\OverGam(D))=\Wrfi(D)$. 
Making the choice $f_i(d)=x_i^d$ we get a map $\Wrxit:{\cal C}\longrightarrow\sympol$
and a commutative diagram with morphisms of graded algebras:


\begin{center}
\setlength{\unitlength}{1pt}
\begin{equation}\label{diagt}
\begin{picture}(150,55)(0,15)
\put(0,50){\makebox(0,0){$A$}}
\put(150,50){\makebox(0,0){${\cal C}$}}
\put(75,-19){\makebox(0,0){$\sympol$}}
\put(0,44){\vector(1,-1){54}}
\put(148,44){\vector(-1,-1){54}}
\put(7,49){\vector(1,0){133}}
\put(75,57){\makebox(0,0){$\OverGam$}}
\put(-7,15){\makebox{$\Wrxi$}}
\put(132,15){\makebox{$\Wrxit$}}
\end{picture}\nopagebreak
\end{equation}\nopagebreak\vspace*{40pt}\nopagebreak
\end{center}\nopagebreak

Now we prove that the weight systems $\Wrxi$ for all natural numbers $r$
are as strong as the universal chromatic weight system.

\begin{lemma}\label{l:iso}
The map $\Wrxit:{\cal C}_{\leq r} \longrightarrow \sympol_{\leq r}$
is an isomorphism of \vs s.
\end{lemma}
{\bf Proof:}
By the definition of $\Wrxit$ we have for single vertices of weight $d$:

\begin{equation}\label{valueofvertex}
\Wrxit \left(
\unitlength=1ex
\begin{picture}(2.0,2.5)(0.0,0.5)
\put(0.7,-0.05){\circle*{0.5}}
\put(0.3,0.95){\makebox(1,1){$d$}}
\end{picture}
\right)=x_1^d+\ldots + x_r^d=G_{d,r}.
\end{equation}

This formula shows that the graded map $\Wrxit$ is surjective and we have

$$
\dim {\cal C}_{\leq r}\leq \dim \sympol_{\leq r}
$$

because $\deg G_{d,r}=d$ and ${\cal C}_{\leq r}$ is spanned by graphs of weight $\leq r$
consisting of isolated vertices.
$\Box$

\medskip

{\bf Proof of Theorem~\ref{linearcomb}}
By Proposition~\ref{WLVr} and Lemma~\ref{coefffkt} the weight system $\Wdxi$ is a 
linear combination of the form described in Theorem~\ref{linearcomb}. 
By the commutative triangle~(\ref{diagt}) 
and Lemma~\ref{l:iso} we see that every weight system of degree $d$ that factors 
through ${\cal C}_d$ also factors through $\sympold_d$. This proves 
Theorem~\ref{linearcomb}.$\Box$

\medskip

As mentioned before the weight system $\Wrfi$ is not a far-reaching generalization
of $\Wrxi$.

\begin{remark}
The weight system $\Wrfi$ is also a chromatic weight system; so for degree $d$,
it factors through $\sympold_d$.
\end{remark}

The restriction to $A_d$ of
the coefficients of monomials in
$x_1,\ldots,x_r$ in $\Wrxi$ are weight systems lying in
the algebra generated by the weight systems of finite degree $W\vert A_i$ 
($i=1,2,\ldots$). Now Lemma~\ref{l:iso} implies immediately the following proposition
that can also be proved more directly.
 
\begin{prop}\label{dpolyiso}
There exists an isomorphism of graded algebras

$$
R[x_1,x_2,\ldots]\longrightarrow \bigoplus_{i=0}^\infty {\cal C}_i^*\ \ (\deg x_i=i)
$$

such that the composition with the inclusion of algebras

$$
\bigoplus_{i=0}^\infty (\OverGam\vert A_i)^*:
\bigoplus_{i=0}^\infty {\cal C}_i^*\longrightarrow \bigoplus_{i=0}^\infty A_i^*
$$

is given by $x_i\mapsto (W\vert A_i)\in A_i^*$.
\end{prop}

The system of generators for the algebra spanned by chromatic weight 
systems of finite degree we used in the proposition consists for $d\geq 2$
not of primitive elements of $\bigoplus_i A_i^*$.

\section{Embeddings into a 2-sphere}\label{embed}

By Section~6.3 of \cite{BN1} there is a connection between the weight systems
$W_{\gl}$ and $W_{\so}$ and surfaces. The description of $W_{\gl}$ and $W_{\so}$
for trivalent graphs instead of chord diagrams allows easily also to consider the 
adjoint representation instead of the defining representation.
We will use \cite{BN2} to give a formula for the value of the universal 
chromatic weight system on certain trivalent diagrams in terms of embeddings into 
a $2$-sphere.

A special case of Proposition~\ref{colorsum1} is the following:

\begin{remark}\label{smiling}
Denote the map $\Wrfi$ for the choice $r=2$ and $f_1=f_2=1$ by $\Wto$. Then the 
formula $\Wto(D)=\# C_2(\Gamma(D))$ for chord diagrams $D$ from 
Proposition~\ref{colorsum1} can be restated as

$$
\Wto (D)=\left\{
\begin{array}{ll}
0 & \mbox{if $\Gamma(D)$ is not bipartite,}\\
2^{\#\pi_0(\Gamma(D))} & \mbox{otherwise.}
\end{array}
\right.
$$
\end{remark}

Let $D$ be a trivalent diagram. We choose a cyclic order for the trivalent vertices 
on $S^1$ such that it is counterclockwise in our pictures. 
If there exists an embedding of $D$
into an oriented $2$-sphere, then let $\sigma(D)$ be the number of trivalent vertices of
$D$ for which the cyclic order does not coincide with the counterclockwise cyclic order
induced by the 
chosen embedding. Let $\sigma(D)$ be arbitrary if no embedding exists. 
Then by \cite{BN2} the formula

\begin{equation}
E(D):=(-1)^{\sigma(D)}\#
\{\mbox{embeddings of $D$ into an oriented $2$-sphere}\}
\end{equation}

defines a weight system.
By Remark~\ref{smiling} we have 

\begin{equation}\label{W21E}
\Wto(D)=E(D)
\end{equation}

for all chord diagrams $D$ hence for all elements of $A$. As an example 
the embeddings of a 
chord diagram $D$ with $E(D)=4$ are shown in Picture~(\ref{edisfour}).

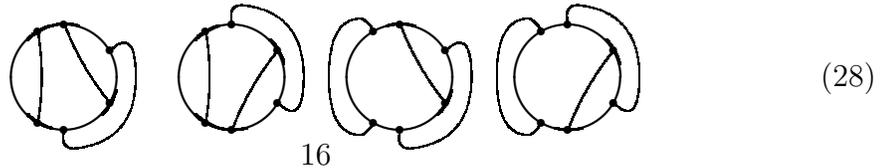
\begin{center}
\begin{equation}\label{edisfour}
\Picture{
\FullCircle
\FullChord[4,8]
\FullChord[3,11]
\ExtChord[1,9]
}
\qquad
\Picture{
\FullCircle
\FullChord[4,8]
\FullChord[9,1]
\ExtChord[3,11]
}
\qquad
\Picture{
\FullCircle
\ExtChord[4,8]
\FullChord[3,11]
\ExtChord[1,9]
}
\qquad
\Picture{
\FullCircle
\ExtChord[4,8]
\FullChord[9,1]
\ExtChord[3,11]
}
\end{equation}\nopagebreak\vspace*{16pt}
\end{center}

Now we can state the 'geometric version' of Theorem~\ref{linearcomb}.

\begin{prop}\label{embeddings}
For $a\in P(A)_d$ we have  

$$
\OverGam(a)=\frac{1}{2}E(a)\,
\unitlength=1ex
\begin{picture}(2.0,2.5)(0.0,0.5)
\put(0.7,-0.05){\circle*{0.5}}
\put(0.3,0.95){\makebox(1,1){$d$}}
\end{picture}.
$$

where $\stackrel{\textstyle d}{\scriptstyle \bullet}$ denotes a single vertex of weight $d$.
\end{prop}
{\bf Proof:}
By definition of $\Wrfi$ we have for $a\in P(A)_d$:

$$
\Wrfi(a)=\sum_{i=1}^r f_i(d)\prod_{{j=1} \atop {j\not=i}}^r f_j(0) W(a).
$$

The previous formula implies $\Wto(a)=2W(a)$ and by Formula~(\ref{W21E}) we have

$$
\Wrfi(a)=\frac{1}{2}\sum_{i=1}^r f_i(d)\prod_{{j=1} \atop {j\not=i}}^r f_j(0) E(a).
$$

In particular, $\Wrxi(a)=\frac{1}{2}(x_1^d+\ldots+x_r^d) E(a)$.
This implies by the commutative 
triangle~(\ref{diagt}), Lemma~\ref{l:iso} and Formula~(\ref{valueofvertex}) 
the desired formula for $\OverGam(a)$.
$\Box$

\medskip

Notice that Proposition~\ref{embeddings} gives us a formula for $\OverGam$ for all 
elements of $\bigcup_{i=0}^\infty {\cal M}^i$ (see Formula~(\ref{primel})) in terms
of $E$.
As an example consider the set

\begin{equation}
{\cal T}:=
\left\{ 
{\mbox{\tt planar trivalent diagrams $D$} \atop 
\mbox{\tt such that $D\setminus S^1$ is a tree}}
\right\}.
\end{equation}

We have ${\cal T}\subset{\cal M}$ 
and elements of ${\cal T}$ have exactly two embeddings
into a 2-sphere.
By Proposition~\ref{embeddings} we have

\begin{equation}
\OverGam(D)=
\unitlength=1ex
\begin{picture}(2.0,2.5)(0.0,0.5)
\put(0.7,-0.05){\circle*{0.5}}
\put(0.3,0.95){\makebox(1,1){$d$}}
\end{picture}.
\end{equation}

for elements $D\in {\cal T}$ of degree $d$.
The element of ${\cal T}$ shown in Picture~(\ref{tree}) is 
for degree $d$ equal to the
element $(-1)^d p_d$ of Proposition~2 of \cite{CDL3}.

\begin{center}
\begin{equation}\label{tree}
\Picture{
\DottedCircle
\Arc[1]\Arc[2]\Arc[3]\Arc[4]\Arc[5]
\FullChord[0,6]\Arc[7]\Arc[8]\Arc[10]\Arc[11]
\Endpoint[7]\Endpoint[8]\Endpoint[10]\Endpoint[11]
\thinlines
\put(-0.866,-0.5){\line(1,2){0.25}}
\put(-0.5,-0.866){\line(1,3){0.289}}
\put(0.866,-0.5){\line(-1,2){0.25}}
\put(0.5,-0.866){\line(-1,3){0.289}}
\put(-1,-1){\makebox(2,1){$\scriptstyle\ldots$}}
}\nopagebreak\vspace*{12pt}
\end{equation}
\end{center}

This leads us to the following remark.

\begin{remark}\label{forest}
The subalgebra of $A$ generated by ${\cal T}$ is also generated by 
the chord diagrams $D$ such that $\Gamma(D)$ is a tree.
It is called the forest algebra in \cite{CDL3}.
\end{remark}

Remark~\ref{forest} follows from \cite{Vo2}, Application on page~7 and 
\cite{CDL3}, Theorem~1 and Proposition~2.

\section{Exceptional Lie superalgebras}\label{other}

In \cite{Vo2} the independence of the weight systems associated to the Lie
superalgebras of Part~III of Table~(\ref{tlsac}) (Theorem 7.1 of \cite{Vo2})
and the insufficience of semisimple Lie superalgebras with Casimir
to generate all weight systems (Theorem 7.4 of \cite{Vo2}) are proved.

\begin{equation}\label{tlsac}
\begin{array}{l}
III \left\{\begin{array}{l}
\\ \\ \\ \\
\end{array}\right.\\
\\ 
\end{array}
\begin{array}{l}
\spli(m,n),\quad (m\not= n)\\
\osp(m,n), \quad \mbox{($n$ even)}\\
\mbox{exceptional Lie algebras $E_6, E_7, E_8, F_4, G_2$}\\
D(2,1,\alpha), \quad (\alpha\in R\setminus\{0,-1\})\\
\mbox{exceptional Lie superalgebras $G(3), F(4)$}\\
\end{array}
\begin{array}{l}
\left.\begin{array}{l}
\\ \\ 
\end{array}\right\} I\\
\left.\begin{array}{l}
\\ \\ \\ 
\end{array}\right\} II
\end{array}
\end{equation}

In Theorem~\ref{linearcomb} we discovered a relation between chromatic weight
systems and the Lie superalgebras of Part~I of
Table~(\ref{tlsac})
\footnote{The weight systems associated to the
standard representation of a Lie superalgebra from Part~I
of Table~(\ref{tlsac}) only depends on the superdimension $m-n$ of the
standard representation, so they are
already determined by $\spli_n$ and $\so_n$.
Using this it can easily be shown that an equivalent formulation of
Theorem~\ref{linearcomb} could use each of the four classical series of
simple Lie algebras.}.
In this section we consider Part~II of the list.
Using the methods from \cite{Vo2} and Proposition~\ref{embeddings} we can
prove the following:

\begin{prop}\label{simpleLie}
Assume that the ground field $R$ is algebraically closed and of
characteristic $0$. For all $d\geq 8$
there exists an element $P_d\in P(A)_d$ with $\OverGam(P_d)\not=0$
and the following property:
Let $\L$ be a simple Lie superalgebra from Part~II of 
Table~(\ref{tlsac}) and
let $\omega$ be a Casimir element for $\L$. Then for all 
finite-dimensional representations $\rho$ of $\L$ we have
$W_{\L,\omega,\rho}(P_d)~=~0$.
\end{prop}
{\bf Proof:}
In \cite{Vo2} a commutative graded $R$-algebra $\Lambda$ and the structure
of a graded $\Lambda$-module on $P(A)_{\geq 2}$ are defined.
Let $\L$ be an exceptional Lie algebra or a Lie superalgebra of type $G(3)$
or $F(4)$
and let $\omega=\sum_i x_i\otimes y_i$ be the Casimir element of $\L$
associated to the Killing form.
Then the endomorphism $\sum_i \ad_\L(x_i)\circ\ad_\L(y_i)$ is equal to
$\id_\L$.
Theorem 6.1 of \cite{Vo2} and the subsequent remarks say that there
exists a unique homomorphism of graded algebras
$\chi_\L:\Lambda\longrightarrow R[t]$
such that for all $\alpha \in\Lambda$, for all $\L$-representations $\rho$
and for all $a\in P(A)_{\geq 2}$ we have

$$
W_{\L,\omega,\rho}(\alpha a)=\chi_\L(\alpha)(1/2)W_{\L,\omega,\rho}(a).
$$

Certain elements $x_1,x_3,x_5,\ldots$ of $\Lambda$ generate a subalgebra $\Lambda_0$, 
on which formulas for the maps $\chi_\L$ are known.
When we define $t:=x_1 /2\in\Lambda_0$ we always have $\chi_\L(t)=t\in R[t]$.
A formula for $\chi_\L(x_n)$ for exceptional Lie algebras is given
by Theorem~6.9 of \cite{Vo2}.
Following Section~6.8 and the remark on page~40 of \cite{Vo2} we
determine for all $n>0$:

$$
\chi_{G(3)}(x_n)=\frac{1}{6}\left(3+2^{n+1}+5(-1)^{n+1}\right)t^n 
$$

and

$$
\chi_{F(4)}(x_n)=\frac{1}{16}\left(8+2^{n+1}+8+27\left(-2/3\right)^{n+1}
\right)t^n.
$$ 

In particular, we have

$$
\chi_{G(3)}(x_3-4t^3)=0 \ \ \mbox{and}\ \ \chi_{F(4)}(3x_3-7t^3)=0.
$$

Let $p_2$ be the trivalent diagram of degree $2$ shown in Picture~(\ref{tree}). 
As in the proof of Theorems~7.1 and 7.4 of \cite{Vo2} we see that for all 
$d\geq 23$ the element

\begin{eqnarray*}
\overtilde{P}_d & = & t^{d-23}(8x_3-7t^3)(81x_3-64t^3)(15x_3-11t^3)(81x_3-85t^3)\\
& & \times (24x_3-41t^3)(x_3-4t^3)(3x_3-7t^3)\cdot p_2
\end{eqnarray*}

has the property $W_{\L,\omega,\rho}(\overtilde{P}_d)=0$ for all simple Lie superalgebras 
we consider in the proposition. The argument concerning the Lie superalgebras of 
type $D(2,1,\alpha)$ is contained in the proof of Theorem~7.4 of \cite{Vo2}.

We can prove that for all monomials $x$ in $t,x_3,x_5 \ldots$ the trivalent
diagram representing $x\cdot p_2$ has exactly two embeddings into a
$2$-sphere.
As examples, for
all degrees $d\geq2$ the element of Picture~(\ref{tree}) is equal to
$t^{d-2}p_2$ and $t^3x_3x_5^2x_7\cdot p_2$ is shown
in Picture~(\ref{bigdiag}).

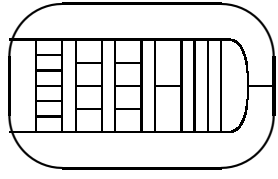
\begin{center}
{\unitlength=25pt
\begin{equation}\label{bigdiag}
\begin{picture}(5,1.4)(-2.5, -0.08)
\thicklines
\put(0,0){\oval(4,2.5)}
\thinlines
\put(-2,-0.7){\line(1,0){3.3}}
\put(-2,0.7){\line(1,0){3.3}}
\put(2,0){\line(-1,0){0.4}}
\qbezier(1.3,0.7)(1.6,0.7)(1.6,0.0)
\qbezier(1.3,-0.7)(1.6,-0.7)(1.6,0.0)
\multiput(1.2,-0.7)(-0.2,0){3}{\line(0,1){1.4}}
\multiput(0.6,-0.7)(-0.4,0){2}{\line(0,1){1.4}}
\put(0.6,0.0){\line(-1,0){0.4}}
\multiput(0,-0.7)(-0.4,0){2}{\line(0,1){1.4}}
\multiput(0,-0.35)(0,0.35){3}{\line(-1,0){0.4}}
\multiput(-0.6,-0.7)(-0.4,0){2}{\line(0,1){1.4}}
\multiput(-0.6,-0.35)(0,0.35){3}{\line(-1,0){0.4}}
\multiput(-1.2,-0.7)(-0.4,0){2}{\line(0,1){1.4}}
\multiput(-1.2,-0.466)(0,0.233){5}{\line(-1,0){0.4}}
\end{picture}
\end{equation}\nopagebreak\vspace*{21pt}
}
\end{center}

This implies by Proposition~\ref{embeddings}:

$$
\OverGam(\overtilde{P}_d)=(8-7)(81-64)(15-11)(81-85)(24-41)(1-4)(3-7)
\unitlength=1ex
\begin{picture}(2.0,2.5)(0.0,0.5)
\put(0.7,-0.05){\circle*{0.5}}
\put(0.3,0.95){\makebox(1,1){$d$}}
\end{picture}
\not=0.
$$ 

and completes the proof for all $d\geq 23$.
The better bound $d\geq 8$ is obtained
by a slightly longer computation using the following sequence of elements:

$$
P_d:=t^{d-8}\cdot (32t^6-71t^3x_3+18x_3^2-45tx_5)\cdot p_2, \quad(d\geq 8).
$$
$\Box$

\medskip

With some technical arguments stated in the Appendix
Proposition~\ref{simpleLie} can be generalized to direct sums of
Lie superalgebras over fields of characteristic $0$.
This will prove the following
complementary result to Theorem~\ref{linearcomb}.

\begin{theorem}
For all degrees $d\geq 8$
the map $\OverGam\vert P(A)_d$ is not a linear 
combination of weight systems coming from
finite direct sums of Lie superalgebras from Part~II of
Table~(\ref{tlsac}).
\end{theorem}

\section{The bialgebra of framed knots}\label{algfknots}

Let $R[\Z/(2)]$ be the group algebra of the group with two elements considered as a
graded bialgebra over $R$ with all elements having degree $0$. Define 

\begin{equation}
\hatA:=A\otimes R[\Z/(2)].
\end{equation}

We can think of elements from $\hatA$ as linear combinations
of chord diagrams, each chord diagram being colored with an element of $\Z/(2)$. We call
the element of $\Z/(2)$ the residue of the diagram.  

Let us explain the origin of the graded bialgebra structure on $\hatA$. By framed
oriented knots we will always mean tame framed oriented knots considered up to
orientation preserving diffeomorphisms of $\R^3$.
Denote the \vs{} freely generated by framed oriented knots by $\fknots$.
The set of framed knots is a semi-group with the connected sum of framed knots as 
multiplication. This turns $\fknots$ into a bialgebra for which framed knots are
group-like elements.
 
Singular framed knots will be considered here as a short
notation for certain linear combinations of framed knots given by the 
'desingularisation rule' from Picture~(\ref{desing}). This rule is valid for three 
singular framed knots that differ only inside a small ball like shown in 
Picture~(\ref{desing}). The framing of the 
shown parts points to the reader.

\begin{center}
\begin{equation}\label{desing}
\Picture{
\thicklines
\put(1,1){\vector(-1,-1){2}}
\put(-1,1){\vector(1,-1){2}}}
\quad = \quad
\Picture{
\thicklines
\put(1,1){\vector(-1,-1){2}}
\put(-1,1){\line(1,-1){0.83}}
\put(0.17,-0.17){\vector(1,-1){0.83}}}
\quad - \quad
\Picture{
\thicklines
\put(1,1){\line(-1,-1){0.83}}
\put(-0.17,-0.17){\vector(-1,-1){0.83}}
\put(-1,1){\vector(1,-1){2}}}
\end{equation}\\[15pt]
\end{center}

Let $\fknots_i$ be the subspace of $\fknots$ spanned by all singular framed 
knots with $i$ double points. This defines a decreasing sequence of \vs s:

\begin{equation}
\fknots=\fknots_0\supseteq \fknots_1\supseteq \fknots_2 
\supseteq\ldots.
\end{equation}

A translation of Proposition~4 from \cite{Vo1} to framed knots tells us that this 
sequence 
is a filtration of the bialgebra $\fknots$. The connected sum of singular framed knots is 
expressed in the same way as for framed knots without double points and the coproduct of 
a singular framed knot can be expressed like in Formula~(3) of \cite{Vo1}.
Denote the graded bialgebra associated to $\fknots$ by 

\begin{equation}
\gr\fknots:=\bigoplus_{i=0}^\infty \fknots_i/\fknots_{i+1}.
\end{equation}

Recall the construction of a singular framed knot $K_D$ as a 'realization' of a chord 
diagram $D$ with residue from \cite{KaT}. We summarize it in Picture~(\ref{realization}).

\begin{center}
\begin{equation}\label{realization}
\mbox{Replace}\quad
\Picture{
\DottedCircle\FullChord[6,0]
\Arc[1]\Arc[2]\Arc[5]\Arc[7]\Arc[8]\Arc[11]}
\quad\mbox{by}\quad\quad
\Picture{
\Arc[1]
\Arc[2]
\Arc[5]
\Arc[6]
\Arc[7]
\Arc[8]
\Arc[11]
\thinlines
\DottedArc[3]
\DottedArc[4]
\DottedArc[9]
\DottedArc[10]
\thicklines
\put(0.9,-0.1){\line(-1,0){1.1}}
\put(-0.2,-0.1){\line(-1,0){0.95}}
\put(0.9,0.1){\line(-1,0){1.725}}
\put(0.707,0.707){\vector(-1,1){0.1}}
\qbezier(0.9,-0.1)(1.04,-0.1)(0.966,-0.259)
\qbezier(0.9,0.1)(1.04,0.1)(0.966,0.259)
\qbezier(-1.15,-0.1)(-1.4,-0.1)(-1.5,0)
\qbezier(-1.15,0.1)(-1.4,0.1)(-1.5,0)
}
\quad \mbox{and} \quad
\Picture{
\DottedCircle\Arc[10]\Arc[11]\Arc[0]\Arc[1]\Arc[2]\Arc[3]
\put(1.1,-0.1){\makebox(0.2,0.2){$\scriptstyle\overline{1}$}}
}
\quad \ \mbox{by} \quad
\Picture{
\Arc[2]\Arc[3]
\Arc[11]\Arc[1]\Arc[10]
\thinlines
\DottedArc[4]\DottedArc[5]\DottedArc[6]\DottedArc[7]
\DottedArc[8]\DottedArc[9]
\thicklines
\qbezier(0.966,0.25)(1.113,-0.25)(1.25,-0.25)
\qbezier(1.25,-0.25)(1.45,-0.25)(1.45,0)
\qbezier(1.25,0.25)(1.45,0.25)(1.45,0)
\qbezier(1.25,0.25)(1.18,0.25)(1.12,0.17)
\qbezier(0.966,-0.259)(0.975,-0.175)(0.98,-0.15)
\put(0.707,0.707){\vector(-1,1){0.1}}
}
\end{equation}\nopagebreak\vspace*{12pt}
\end{center}

If $\deg D=d$, then $K_D\in \fknots_d/\fknots_{d+1}$ is well defined.
This gives us an isomorphism of graded bialgebras
$\realize:\hatA \longrightarrow \gr\fknots.$
The special case $\dim \fknots_0/\fknots_1=2$ may be
proved directly. In general 
the proof that $\realize$ is injective is based on Part~(1) of the next theorem.
In order to formulate it let us introduce some notation.
We can pass from a graded \vs{} $B=\bigoplus_{i=0}^\infty B_i$ to its completion 
$B^c:=\prod_{i=0}^\infty B_i$. We will use the notation $\sum_{i=0}^\infty b_i h^i$ 
with $b_i\in B_i$ and a formal parameter $h$ for elements of $B^c$. If $B$ is a graded 
algebra, then $B^c$ is an algebra.
The elements of $B^c$ are multiplied like formal power series.

\begin{theorem}[Kontsevich]\label{Kontsevich}
(1) There exists a linear map
$Z:\fknots \longrightarrow \hatAc$

such that for any realisation $K_D$ of a chord diagram $D$ with residue we have

$$
Z(K_D)=D h^{\deg D}+\mbox{terms of higher degree}.
$$

(2) With the comultiplication defined by 
$\Delta(\sum_{i=0}^\infty a_i h^i):=\sum_{i=0}^\infty \Delta(a_i) h^i$ the image of the
map $Z$ is a bialgebra and $Z:\fknots\longrightarrow \image Z$
is a morphism of bialgebras.
\end{theorem}

For a proof of Part~(1) of Theorem \ref{Kontsevich} see \cite{KaT} and also \cite{BN1}
and \cite{LM2}. For a proof of Part~(2) and a proof that $Z$ is already defined over 
the rational numbers one may use the explicit description of $Z$
from \cite{KaT} and Theorem $A''$ from \cite{Dri}.
We have used a normalization of $Z$ different from the one used in \cite{KaT} 
in order to turn $Z$ into a morphism of algebras\footnote{We give more details in the 
proof of Proposition~\ref{Lymphtofu}.}. Our
normalization specializes to the one in \cite{BN1} when we forget the residues
and introduce the 'framing independence' relation (see \cite{LM2}, Theorem~6 
and Theorem~7).
For reasons that will become clear in the next section the map $Z$ is
called the {\em universal Vassiliev invariant}.
It is an open question whether $Z$ 
is injective.

\section{Vassiliev invariants}\label{vi}

Now we take a look at the dual notions of the filtered bialgebra of framed knots and
the graded bialgebra of chord diagrams and translate the results of the previous section.

\begin{defi}
We call a linear map $v:\fknots\longrightarrow R$ a Vassiliev invariant of degree $d$
if $v(\fknots_{d+1})=\{0\}$. Let ${\cal V}_{d}$ be the space of Vassiliev invariants
of degree $d$ and ${\cal V}:=\bigcup_{i=0}^\infty {\cal V}_d$ the space of all Vassiliev
invariants.
\end{defi}

The \vs s ${\cal V}_i$ form an increasing sequence

\begin{equation}
{\cal V}_0\subseteq {\cal V}_1\subseteq {\cal V}_2\subseteq 
\ldots \subseteq {\cal V}.
\end{equation}

The product of $v_1\in{\cal V}_k$ and $v_2\in{\cal V}_l$ is defined to be pointwise on 
framed knots $K$ without double points: $(v_1 v_2)(K):=v_1(K)v_2(K).$
With this definition ${\cal V}$ becomes a subalgebra of ${\fknots}^*$. Now
$v_1 v_2$ is in ${\cal V}_{k+l}$ because the coproduct of $\gr \fknots$ respects the 
grading.
Define ${\cal V}_{-1}:=\{ 0 \}$. 
Theorem \ref{Kontsevich} enables us to choose nice representatives for elements of 
${\cal V}_d/{\cal V}_{d-1}$ ($d\geq 0$). This will allow us to turn ${\cal V}$ into a
graded algebra as follows:
Denote the projection $\hatAc\longrightarrow \hatA_d$ by $p_d$. Given $w\in \hatA_d^*$, we 
can define $\VI_d(w) \in {\cal V}_d$ by

\begin{equation}
\VI_d(w)(K)=(w\circ p_d\circ Z)(K).
\end{equation}

The image $\Vkan_d:=\VI_d(\hatA_d^*)$ of the linear map $\VI_d$ is called the
space of canonical Vassiliev invariants of degree $d$.
We state the translation of Theorem~\ref{Kontsevich} to this setting 
(see \cite{BN1}, Theorem~9).

\begin{prop}\label{kaniso}
(1) The map

$$
\bigoplus_{i=0}^\infty \VI_i:\bigoplus_{i=0}^\infty \hatA_i^*\longrightarrow
\bigoplus_{i=0}^\infty \Vkan_i={\cal V}
$$

is an isomorphism of algebras.

(2) The image of the space of primitive elements 
$P(\bigoplus_{i=0}^\infty \hatA_i^*)$ is the space of Vassiliev invariants that are 
additive with respect to the connected sum of framed oriented knots.
\end{prop}

Now we pass to completions.

\begin{defi}
A formal sum $v=\sum_{i=0}^\infty v_ih^i$ with $v_i\in {\cal V}_i$ is called a Vassiliev
series. The Vassiliev series $v$ is called canonical if 
$v\in\prod_{i=0}^\infty \Vkan_i$.
The Vassiliev series $v$ is called multiplicative if the map 
$v:\fknots\longrightarrow R[[h]]$ is a morphism of algebras.
\end{defi}

We have a projection $\forget_d:\hatA_d\longrightarrow A_d$ defined by forgetting residues.
The dual map allows to lift elements from $A_d^*$ to $\hatA_d^*$.

\begin{defi}
(1) For a weight system $w$ define a canonical Vassiliev series by

$$
\VS(w):=\sum_{i=0}^\infty (w\circ\forget_i\circ p_i\circ Z)h^i
$$

(2) For a Vassiliev series $v=\sum_{i=0}^\infty v_i h^i$ define a weight system 
$\WS(v)\in A^*$ by $\WS(v)(D):=v_d(K_D)$ where $D$ is a chord diagram of degree $d$ and
$K_D$ is a realization of $D$ with residue $\overline{0}$.
\end{defi}

By Part~(1) of Theorem~\ref{Kontsevich} we have 
for every weight system $w$ : 

\begin{equation}\label{psiphi}
\WS(\VS(w))=w. 
\end{equation}

A consequence of Part~(2) of 
Theorem~\ref{Kontsevich} is the following proposition (see also \cite{BNG}, 
Proposition~2.9).

\begin{prop}\label{p:VSandWS}
(1) The canonical Vassiliev series $\VS(w)$ is multiplicative if
and only if $w$ is multiplicative. 

(2) For weight systems $w_1$ and $w_2$ we 
have $\VS(w_1 w_2)=\VS(w_1)\VS(w_2)$.
\end{prop}

\section{The knot polynomial $\V1$}

Define a linear map $\V1: \fknots\longrightarrow R[y,y^{-1}]$ by:

\begin{eqnarray*}
& (1) & \V1(K_+)-\V1(K_-)=(y-y^{-1})\V1(K_{\mid\mid 1})\V1(K_{\mid\mid 2}),\\
& (2) & \V1(K_t)=y\V1(K_\mid),\\
& (3)  & \V1(O)=1
\end{eqnarray*}

where we use two knots $K_+$ and $K_-$ and a link $K_{\mid\mid}$ with
two connected components that differ only in a ball like shown in 
Picture (\ref{homkauff}) for the first equation. 
We denote the two knots that are the connected components
of $K_{\mid\mid}$ by $K_{\mid\mid 1}$ and $K_{\mid\mid 2}$. The symbols $K_t$ and
$K_\mid$ from the
second equation are also explained by Picture (\ref{homkauff}). The symbol $O$ from
the third equation denotes the oriented unknot with $0$-framing. 
Proposition~\ref{VSWandV1} says that $\V1$ is well defined. 

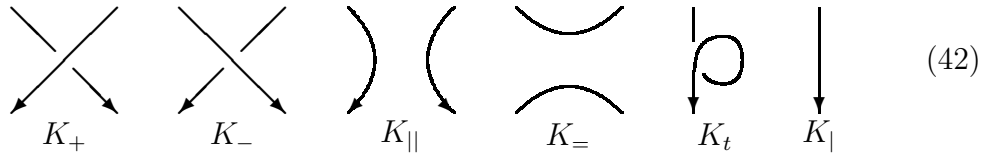
\begin{center}
\begin{center}
\begin{equation}\label{homkauff}
\Picture{
\thicklines
\put(1,1){\vector(-1,-1){2}}
\put(-1,1){\line(1,-1){0.83}}
\put(0.17,-0.17){\vector(1,-1){0.83}}
\put(-1,-1.8){\makebox(2,0.7){$K_+$}}
}
\qquad
\Picture{
\thicklines
\put(1,1){\line(-1,-1){0.83}}
\put(-0.17,-0.17){\vector(-1,-1){0.83}}
\put(-1,1){\vector(1,-1){2}}
\put(-1,-1.8){\makebox(2,0.7){$K_-$}}
}
\qquad
\Picture{
\thicklines
\qbezier[90](1,1)(0,0)(1,-1)
\qbezier[90](-1,1)(0,0)(-1,-1)
\put(-0.9,-0.9){\vector(-1,-1){0.1}}
\put(0.9,-0.9){\vector(1,-1){0.1}}
\put(-1,-1.8){\makebox(2,0.7){$K_{\mid\mid}$}}
}
\qquad
\Picture{
\thicklines
\qbezier[90](1,1)(0,0)(-1,1)
\qbezier[90](1,-1)(0,0)(-1,-1)
\put(-1,-1.8){\makebox(2,0.7){$K_=$}}
}
\qquad
\begin{picture}(1.2,1)(-1,-0.167)
\thicklines
\put(-0.8,1){\line(0,-1){0.64}}
\put(-0.8,-0.2){\vector(0,-1){0.8}}
\qbezier[50](-0.63,-0.27)(-0.53,-0.45)(-0.25,-0.45)
\qbezier[50](-0.25,-0.45)(0.1,-0.45)(0.1,0.0)
\qbezier[50](-0.25,0.45)(0.1,0.45)(0.1,0.0)
\qbezier[50](-0.8,-0.2)(-0.8,0.45)(-0.25,0.45)
\put(-1,-1.8){\makebox(1.2,0.7){$K_t$}}
\end{picture}
\qquad
\begin{picture}(0.4,1)(-0.2,-0.167)
\thicklines
\put(0,1){\vector(0,-1){2}}
\put(-0.5,-1.8){\makebox(1,0.7){$K_\mid$}}
\end{picture}
\end{equation}\nopagebreak\vspace*{21pt}
\end{center}
\end{center}

Define $\V1_h(K):=\V1(K)(\exp(h/2))$ (this shall mean that $y^{\pm 1}$ is
substituted by $\exp(\pm h/2)$). Then $\V1_h$ is a power series and the 
coefficients of powers of $h$ are knot invariants with values in $R$.
Now we formulate our second main result.

\begin{theorem}\label{chrompol}
The algebra of all chromatic weight systems of finite degree is mapped by $\VS$ to the
subalgebra of $\cal V$ generated by the coefficients of $\V1_h$. 
In particular, the Vassiliev invariants corresponding to all chromatic weight systems 
distinguish exactly the same knots as the polynomial $\V1$. 
\end{theorem}

Theorem~\ref{chrompol} will follow from Proposition~\ref{dpolyiso}, Part~(2)
of Proposition~\ref{p:VSandWS} and Proposition~\ref{VSWandV1}. 
Let us now consider an example.

{
\unitlength=14pt
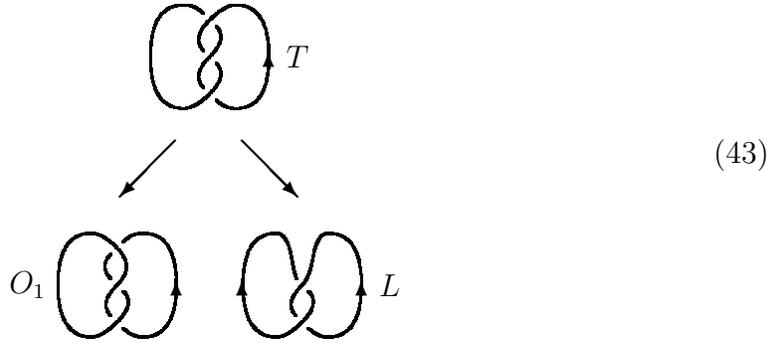
\begin{center}
$$
\begin{picture}(5,2.8)(-2.5,-1.4)
{
\thicklines
\qbezier[50](0,0)(0.5,0.5)(0.17,0.83)
\qbezier[80](-0.17,1.17)(-0.4,1.4)(-0.7,1.4)
\qbezier[80](-0.7,1.4)(-1.6,1.4)(-1.6,0)
\qbezier[80](0,-1)(-0.4,-1.4)(-0.7,-1.4)
\qbezier[80](-0.7,-1.4)(-1.6,-1.4)(-1.6,0)
\qbezier[60](0,-1)(0.5,-0.5)(0.17,-0.17)
\qbezier[60](-0.17,0.17)(-0.5,0.5)(0,1)
\qbezier[80](0.17,-1.17)(0.4,-1.4)(0.7,-1.4)
\qbezier[80](0.7,-1.4)(1.6,-1.4)(1.6,0)
\qbezier[80](0,1)(0.4,1.4)(0.7,1.4)
\qbezier[80](0.7,1.4)(1.6,1.4)(1.6,0)
\qbezier[50](-0.17,-0.83)(-0.5,-0.5)(0,0)
\put(1.62,-0.1){\vector(0,1){0.25}}
\put(1.9,-0.5){\makebox(1,1){$T$}}
}
\end{picture}
$$
\end{center}
\begin{equation}\label{ytree}
\begin{picture}(10,0.2)(-5,-0.1)
\thicklines
\put(-0.9,0.5){\vector(-1,-1){1.5}}
\put(0.9,0.5){\vector(1,-1){1.5}}
\end{picture}
\end{equation}
\begin{center}
$$
\begin{picture}(5,2.8)(-2.5,-1.4)
{
\thicklines
\qbezier[50](0,0)(0.5,0.5)(0.0,1)
\qbezier[80](0,1)(-0.4,1.4)(-0.7,1.4)
\qbezier[80](-0.7,1.4)(-1.6,1.4)(-1.6,0)
\qbezier[80](0,-1)(-0.4,-1.4)(-0.7,-1.4)
\qbezier[80](-0.7,-1.4)(-1.6,-1.4)(-1.6,0)
\qbezier[50](0,-1)(0.5,-0.5)(0.17,-0.17)
\qbezier[50](-0.17,0.17)(-0.5,0.5)(-0.17,0.83)
\qbezier[80](0.17,-1.17)(0.4,-1.4)(0.7,-1.4)
\qbezier[80](0.7,-1.4)(1.6,-1.4)(1.6,0)
\qbezier[80](0.17,1.17)(0.4,1.4)(0.7,1.4)
\qbezier[80](0.7,1.4)(1.6,1.4)(1.6,0)
\qbezier[50](-0.17,-0.83)(-0.5,-0.5)(0,0)
\put(1.62,-0.1){\vector(0,1){0.25}}
\put(-2.9,-0.5){\makebox(1,1){$O_1$}}
}
\end{picture}
\begin{picture}(5,2.8)(-2.5,-1.4)
{
\thicklines
\qbezier[50](0,0)(0.2,0.2)(0.3,0.8)
\qbezier[80](0.3,0.8)(0.42,1.4)(0.7,1.4)
\qbezier[50](-0.15,0.17)(-0.2,0.22)(-0.3,0.8)
\qbezier[80](-0.3,0.8)(-0.42,1.4)(-0.7,1.4)
\qbezier[80](-0.7,1.4)(-1.6,1.4)(-1.6,0)
\qbezier[80](0,-1)(-0.4,-1.4)(-0.7,-1.4)
\qbezier[80](-0.7,-1.4)(-1.6,-1.4)(-1.6,0)
\qbezier[50](0.17,-0.17)(0.5,-0.5)(0,-1)
\qbezier[80](0.17,-1.17)(0.4,-1.4)(0.7,-1.4)
\qbezier[80](0.7,-1.4)(1.6,-1.4)(1.6,0)
\qbezier[80](0.7,1.4)(1.6,1.4)(1.6,0)
\qbezier[50](-0.17,-0.83)(-0.5,-0.5)(0,0)
\put(1.62,-0.1){\vector(0,1){0.25}}
\put(-1.6,-0.1){\vector(0,1){0.25}}
\put(1.9,-0.5){\makebox(1,1){$L$}}
}
\end{picture}
$$
\end{center}
}

Using the tree in Picture~(\ref{ytree}) we can determine the value of $\V1$ on the
shown trefoil knot $T$ as follows. By Relations~(2) and (3) we have $\V1(O_1)=y$. The
two connected components of the link $L$ from Picture~(\ref{ytree}) are unknots $O$ with
$0$-framing, so we use Relations~(1) and (3) to compute 

\begin{equation}\label{V1T}
\V1(T)=\V1(O_1)+(y-y^{-1})\V1(O)\V1(O)=2y-y^{-1}.
\end{equation}

The following remark is easy to prove.

\begin{remark}
Let $K$ be an oriented framed knot with mirror image $\overline{K}$. Then we have

\begin{equation}
\V1(K)(y)=\V1\left(\overline{K}\right)(y^{-1})
\end{equation}
\end{remark}

From the remark and Formula~(\ref{V1T}) we see that the polynomial $\V1$ can 
distinguish at least one framed knot from its mirror image.
\footnote{
The knot invariant $\V1$ has also been studied in Section~4 of \cite{Kne}
in a framing independent normalization. Corollary~2 of \cite{Kne} implies
that the set of all polynomials that appear as values $\V1(K)$ for some  
oriented framed knot $K$ is equal to
$\bigcup_{i\in\Z} y^i(1+(y^2-1)\Z[y^2,y^{-2}])$.}

\section{Canonical Vassiliev series coming from $\gl_n$}\label{Homfly}

Now we translate our results about weight systems into results
about the corresponding canonical Vassiliev series. The next
proposition is proved in \cite{LM1} for links without framing 
with analytical me\-thods. 
We use the same ideas to give a proof for framed links with the help of the algebraic
description of the universal Vassiliev invariant $Z$.
Recall the notation from Picture~(\ref{homkauff}).

\begin{prop}\label{Lymphtofu}
a) For all $n\geq 1$ there exists an invariant $H$ of framed oriented links with values 
in $R[[h]]$ satisfying the following three relations:

\begin{eqnarray*}
& (1) & H(K_+)-H(K_-)=\left(\exp(h/2)-\exp(-h/2)\right) H(K_{\mid\mid}),\\
& (2) & H(K_t)=\exp(nh/2) H(K_{\mid}),\\
& (3) & H(O)=1.
\end{eqnarray*}

b) For framed oriented knots $K$ we have 

$$
H(K)=\VS(W_{\gl_n,X,\tau})(K)/n.
$$ 

c) Define $[n]_q:=\frac{q^n-q^{-n}}{q-q^{-1}}$ for invertible elements $q\in R[[h]]^*$.
Then we have for the disjoint union 
$K_1\amalg K_2$ of framed oriented links $K_1$ and $K_2$:

$$
H(K_1 \amalg K_2)=[n]_{\exp(h/2)} H(K_1)H(K_2).
$$
\end{prop}
{\bf Proof:}
We start with a remark concerning generalisations of Theorem~\ref{Kontsevich} and
Lemma~\ref{welldef} of this article. As a general reference we use
\cite{KaT} and also Chapters~XI and XX of \cite{Kas}. 

The objects of all categories that we consider in this proof are finite sequences of
'$+$'- and '$-$'-symbols.
There exists a category $\hatAtan[[h]]$ in which the morphisms are represented
by formal power series of chord diagrams on circles and intervals with residue 
(see Section~2 of \cite{KaT} for a precise definition of the diagrams we use here and 
for a formula for the residue of the composition of diagrams). 

The universal Vassiliev invariant $\Ztan$ is a functor
from the category ${\cal T}$ of framed oriented tangles to $\hatAtan[[h]]$
In ${\cal T}$ we have links as the 
endomorphisms of the empty sequence and for knots $K$ we have

\begin{equation}\label{normali}
Z(K)=\Ztan(K)/\Ztan(O).
\end{equation}

In Equation~(\ref{normali}) the multiplication is induced by the
connected sum of chord diagrams with residue and we use the inclusion of $\hatA^c$ into
the endomorphisms of the empty sequence in $\hatAtan[[h]]$.

The weight system $W_{\L,\omega, \rho}$ associated to a Lie superalgebra $\L$ with Casimir
element $\omega$ and representation $\rho:\L\longrightarrow \End M$ can be generalized
to a functor also called $W_{\L,\omega, \rho}$ from the category 
$\hatAtan[[h]]$ to a category $D_\L(M,M^*)[[h]]$.
In order to describe the morphisms in the last category we associate to an object
$(\epsilon_1,\ldots,\epsilon_k)$ ($\epsilon_i\in\{+,-\})$ of 
$D_\L(M,M^*)[[h]]$ the module $M^{\epsilon_1}\otimes\ldots\otimes 
M^{\epsilon_k}$ where $M^+$ denotes $M$ and $M^-$ denotes the dual module $M^*$.
The 'empty tensor product' will be $R$.
The morphisms between two objects of $D_\L(M,M^*)[[h]]$ are formal power series of 
$\L$-linear maps between the corresponding modules
In the sequel we assume $\End_\L (M)\cong R$. 
Then by Formula~(\ref{normali}) and Proposition~2.1.6 of \cite{LM1} we have for framed 
knots $K$:

\begin{equation}\label{WZWZtan}
\VS(W_{\L,\omega, \rho})(K)=\frac{\dim M}{W_{\L,\omega, \rho}(\Ztan(O))}
W_{\L,\omega, \rho}(\Ztan(K)).
\end{equation}

We use the canonical isomorphism 
$R[[h]]\cong\End_{D_\L(M,M^*)[[h]]}(\emptyseq)$ in Equation~(\ref{WZWZtan}) where 
$'\emptyseq'$ denotes the empty sequence.

After this long preparation we start with a simple proof of the proposition.
We divide each link $K_x$ ($x\in\{+,-,\mid\mid\}$) from Relation~(1) 
for $H$ into three tangles $T_1$, $T_x$ and $T_2$
as shown in Picture~(\ref{linkproj}).
The tangles $T_1$ and $T_2$ are the same in all three cases and the part of the
tangle $T_x$ denoted by a box labeled $x$ in Picture~(\ref{linkproj}) looks like the 
corresponding part in Picture~(\ref{homkauff}).

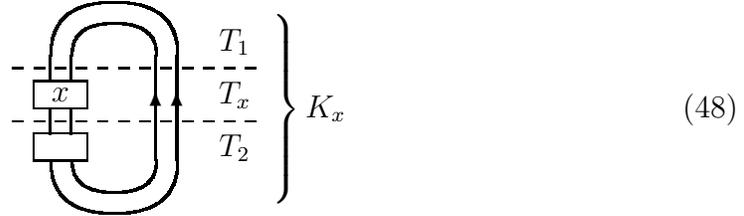
\begin{center}
\begin{equation}\label{linkproj}
\left.
\unitlength=20pt
\begin{picture}(4.8,2)(-1.8,-0.167)
\thicklines
\qbezier(-0.8,1)(-0.8,1.6)(0,1.6)
\qbezier(0.8,1)(0.8,1.6)(0,1.6)
\qbezier(-1.2,1)(-1.2,2)(0,2)
\qbezier(1.2,1)(1.2,2)(0,2)
\put(-1.2,1){\line(0,-1){0.5}}
\put(-0.8,1){\line(0,-1){0.5}}
\thinlines
\put(-1.5,0){\framebox(1,0.5){$x$}}
\put(-1.5,-0.99){\framebox(1,0.5){}}
\thicklines
\put(-1.2,0){\line(0,-1){0.5}}
\put(-0.8,0){\line(0,-1){0.5}}
\put(1.2,0){\line(0,1){1}}
\put(0.8,0){\line(0,1){1}}
\put(1.2,-1){\vector(0,1){1.35}}
\put(0.8,-1){\vector(0,1){1.35}}
\qbezier(-0.8,-1)(-0.8,-1.6)(0,-1.6)
\qbezier(0.8,-1)(0.8,-1.6)(0,-1.6)
\qbezier(-1.2,-1)(-1.2,-2)(0,-2)
\qbezier(1.2,-1)(1.2,-2)(0,-2)
\thinlines
\multiput(-1.9,-0.25)(0.4,0){12}{\line(1,0){0.2}}
\multiput(-1.9,0.75)(0.4,0){12}{\line(1,0){0.2}}
\put(1.8,0.75){\makebox(1,1){$T_1$}}
\put(1.8,-0.25){\makebox(1,1){$T_x$}}
\put(1.8,-1.25){\makebox(1,1){$T_2$}}
\end{picture}
\right\} K_x
\end{equation}\vspace*{24pt}
\end{center}

Now we make the special choice $\L=\gl_n$, $\omega=X$ and $\rho=\tau$ from 
Lemma~\ref{Wgltau}. We define $H$ for framed oriented links $K$ by

\begin{equation}\label{defh}
H(K):=\frac{\Wglnxt(\Ztan(K))}{\Wglnxt(\Ztan(O))}.
\end{equation}

Define $\L$-linear maps $\varphi_i:=\Wglnxt(\Ztan(T_i))$ and

$$
\varphi_x:=\Wglnxt(\Ztan(T_x))\in \End(M^{\otimes 2}\otimes {M^*}^{\otimes 2})[[h]].
$$

By the explicit description of the universal Vassiliev invariant $\Ztan$ (see
Theorem~4.7 of \cite{KaT}) we have

\begin{eqnarray*}
\varphi_+ & = & \left(\exp(hX/2)\circ X\right)\otimes \id_{M^*}\otimes \id_{M^*},\\
\varphi_- & = & \left(\exp(-hX/2)\circ X\right)\otimes \id_{M^*}\otimes \id_{M^*},\\
\varphi_{\mid\mid} & = & \id_M\otimes\id_M\otimes\id_{M^*}\otimes \id_{M^*}.
\end{eqnarray*}

In these formulas no associators appear
because in Picture~(\ref{linkproj}) the box labeled $x$ is on the left side.
A short computation in the algebra $(\End M \otimes \End M)[[h]]$ yields

$$
\left(\exp(hX/2)-\exp(-hX/2)\right) \circ X 
=\left(\exp(h/2)-\exp(-h/2)\right) \id\otimes \id.
$$

This implies Relation~(1) because we have
$(\varphi_1\circ \varphi_x\circ \varphi_2)/\Wglnxt(\Ztan(O))=H(K_x)$.
Property~(2) is proved by a similar computation and Property~(3)
follows from Equation~(\ref{defh}). 

By Equation~(\ref{WZWZtan}) the invariant $H$ is an extension of $\VS(\Wglnxt)/n$
from knots to links. This proves Part b) of the theorem.

A computation for the disjoint union of two unknots 
using Relations (1)--(3) yields
$H(O \amalg O)=[n]_{\exp(h/2)}H(O)^2$.

The functoriality of $\Wglnxt\circ \Ztan$ now allows to determine

\begin{equation}\label{qdim}
\Wglnxt(\Ztan(O))=[n]_{\exp(h/2)}
\end{equation}

and gives us the general formula for the disjoint union of links for $H$
\footnote{By Theorem~7.2 of \cite{KaT} or Theorem~10 of \cite{LM2} we see that the
factor $[n]_{\exp(h/2)}$ from Formula~(\ref{qdim}) appears also as a quantum dimension. 
For the more direct computations usin $R$-matrices see~\cite{Tur}.}.
$\Box$

\medskip

This again proves that the Homfly polynomial 
is well defined. 
The next proposition is a consequence of a connection between $H$ and $\V1$.

\begin{prop}\label{VSWandV1}
The invariant $\V1$ is well-defined and for all framed oriented knots $K$ we have

$$
\VS(W)(K)=\V1(K)(\exp(h/2)).
$$
\end{prop}
{\bf Proof:}
By Lemma \ref{Wgltau} we have for all $n$:

$$
\VS(W_{gl_n,X,\tau})/n=\sum_{i=0}^\infty\sum_{j=0}^i p_{ij}n^j h^i
$$

with uniquely determined Vassiliev invariants $p_{ij}$ of degree $i$ and

\begin{equation}\label{diagsum}
\VS(W)=\sum_{i=0}^\infty p_{ii} h^i.
\end{equation}

From the Relations (1)-(3) for $H$ from Lemma \ref{Lymphtofu} we can derive 
relations for $\VS(W)$ as follows:

Recall from Part~I, Section~5 of \cite{Kau} that we can calculate $H(K)$ for every
link $K$ with the help of a certain labeled binary tree with root, a so-called skein-tree. 
The vertices of such a tree are labeled with link diagrams, the root is labeled
with a diagram for $K$, the leaves are labeled with diagrams of unknotted unlinked
framed circles, and for the link diagrams labeling inner vertices there exists a
crossing such that the relation shown in Picture~(\ref{skeintree}) is satisfied.
An example for a part of a skein-tree is shown in Picture~(\ref{ytree}).

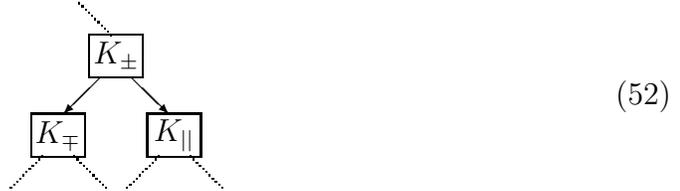
\begin{center}
\begin{equation}\label{skeintree}
\begin{picture}(4,1.8)(-2.1,-0.5)
\thinlines
\put(-0.5,0){\framebox(1,0.8){$K_\pm$}}
\put(-1.6,-1.5){\framebox(1,0.8){$K_\mp$}}
\put(0.6,-1.5){\framebox(1,0.8){$K_{\mid\mid}$}}
\put(-0.3,0){\vector(-1,-1){0.7}}
\put(0.3,0){\vector(1,-1){0.7}}
\qbezier[10](-0.1,0.8)(-0.4,1.1)(-0.7,1.4)
\qbezier[10](-1.4,-1.5)(-1.7,-1.8)(-2.0,-2.1)
\qbezier[10](-0.8,-1.5)(-0.5,-1.8)(-0.2,-2.1)
\qbezier[10](1.4,-1.5)(1.7,-1.8)(2.0,-2.1)
\qbezier[10](0.8,-1.5)(0.5,-1.8)(0.2,-2.1)
\end{picture}
\end{equation}\vspace*{18pt}
\end{center}

Notice that $H(K_{\mid\mid})$ is multiplied in Relation (1) with

$$
\exp(h/2)-\exp(-h/2) = h+\mbox{terms of degree $>1$}
$$

and for the factor $[n]_{\exp(h/2)}$ for the disjoint union of links we have 

$$
\frac{\exp(nh/2)-\exp(-nh/2)}{\exp(h/2)-\exp(-h/2)} = 
n \sum_{\nu=0}^\infty \frac{(nh/2)^{2\nu}}{(2\nu+1)!}
+\mbox{terms $a_{ij} n^ih^j$ with $i\leq j$}.
$$

This implies that in a skein-tree 
for the calculation of $H(K)$ the subtree for the link $K_{\mid\mid}$ can only give a
contribution to the highest powers in $n$ if the number of components of
$K_{\mid\mid}$ is greater than the number of components of $K_+$ (and $K_-$). This is
the case exactly if the two strands of the distinguished crossing of the link $K_+$
belong to the same connected component. 
Using this argument we see that $\VS(W)$ satisfies the Relations (1)-(3) for $\V1$ with 
the parameter $y=\exp(h/2)$.
$\Box$

\medskip

This completes the proof of Theorem~\ref{chrompol}.
From Part~(1) of Proposition~\ref{p:VSandWS} and Formula~(\ref{psiphi}) we get the 
following corollary to Proposition~\ref{VSWandV1}.

\begin{coro}
(1) The linear map $\V1$ is a morphism of algebras.

(2) The power series $\V1_h$ is a canonical Vassiliev series
with corresponding weight system $\WS(\V1_h)=W$.
\end{coro}

A reason why we mainly considered knots and not links is given in the following remark
that can be proved in the same way as Proposition~\ref{VSWandV1}.

\begin{remark}
A generalization of $\V1$ to links suggested by Formula~(\ref{diagsum})
is not useful because we would have $\V1(L)=\V1(K_1)\cdot\ldots\cdot\V1(K_l)$ where
$K_1,\ldots,K_l$ are the knots that are the connected components of the link $L$. 
\end{remark}

Now consider the Lie algebra $\so_n$ with Casimir element $\omega$ and representation
$\tau$ as in Lemma \ref{Wso}. Recall again the symbols 
from Picture~(\ref{homkauff}).

\begin{prop}\label{kauffman} 
a) For all $n\geq 2$ there exists an invariant $\K$ of framed oriented links with values 
in $R[[h]]$ satisfying the following three relations:

\begin{eqnarray*}
& (1) & \K(K_+)-\K(K_-)=\left(\exp(h/2)-\exp(-h/2)\right)
\left(\K(K_{\mid\mid})-\K(K_=)\right),\\
& (2) & \K(K_t)=\exp((n-1)h/2)\K(K_\mid),\\
& (3) & \K(O)=1
\end{eqnarray*}

where we may choose any orientation for $K_=$. 

b) For framed oriented knots $K$ we have 

$$
\K(K)=\VS(W_{\so_n,\omega,\tau})(K)/n.
$$

c) For the disjoint union of framed oriented links $K_1$ and $K_2$ we have

$$
\K(K_1 \amalg K_2)=\left([n-1]_{\exp(h/2)}+1\right) \K(K_1)\K(K_2).
$$
\end{prop}

A proof of Proposition~\ref{kauffman} is only slightly more complicated than the proof
of Proposition~\ref{Lymphtofu} (see Theorem~3.6 of \cite{LM3}).
The invariant $\K$ does not depend on the orientation of framed links. 
It is one of the standard parametrizations of the Kauffman polynomial
(see \cite{Kau}). 
Because of Lemma \ref{Wso} we could also use it instead of $H$ to prove 
Proposition \ref{VSWandV1}. In particular, $\V1$ is also a specialization
of the Kauffman polynomial.

\appendix

\section{Appendix}

In order to generalize Proposition~\ref{simpleLie} from a simple Lie
superalgebra to a direct
sum $\L=\bigoplus_{i=1}^k \L_i$ of simple Lie superalgebras $\L_i$ we have to use a 
relation between the description of $\L$-modules and the description of
$\L_i$-modules.
A technical problem is that a finite-dimensional $\L$-module need not be the direct
sum of simple $\L$-modules. 
To solve this problem recall the definition of a Grothendieck ring 
(see Part~16B of \cite{CuR}). 
By $G_0(\L)$ we will denote the Grothendieck ring associated to the category
of finite dimensional $\L$-modules using all short exact sequences of $\L$-modules
for the defining relations in $G_0(\L)$. For the next lemma see also Remark~2.12 of
\cite{BN1}.

\begin{lemma}
The weight systems associated to a Lie superalgebra $\L$ with a fixed Casimir element
$\omega$ induce a morphism of abelian groups

$$
W_{\L,\omega,\cdot}:G_0(\L)\longrightarrow A^* \ ,\ \ [\rho]\mapsto W_{\L,\omega,\rho}
$$
\end{lemma}
{\bf Proof:}
Let $M'\subset M$ be finite dimensional $\Z/(2)$-graded vector-spaces and let
$f:M\longrightarrow M$ be a
linear map satisfying $f(M')\subseteq M'$. Denote the restriction of $f$ to $M'$ by 
$f'$ and the induced endomorphism of $M'':=M/M'$ by $f''$. 
Then the supertrace satisfies $\str(f')+\str(f'')=\str(f)$. This completes the proof
because $W_{\L,\omega,\rho}(D)$ can be described as the supertrace of $\rho(a_D)$
where $a_D$ is an element of the universal enveloping algebra $U(\L)$ depending only 
on $D$ and $\omega$ but not on $\rho$.
$\Box$

\medskip

Now we can generalize Proposition~\ref{simpleLie} to semisimple Lie superalgebras
over fields of characteristic $0$.

\begin{prop}\label{semisimpleLie}
Let the ground field $R$ be an arbitrary field of characteristic~$0$.
For all $d\geq 8$
there exists an element $P_d\in P(A)_d$ with $\OverGam(P_d)\not=0$
and the following property:
Let $\L$ be a finite direct sum of Lie superalgebras from Part~II of
Table~(\ref{tlsac})
and let $\omega$ be a Casimir element for $\L$. Then for all 
finite-dimensional representations $\rho$ of $\L$ we have $W_{\L,\omega,\rho}(P_d)=0$.
\end{prop}
{\bf Proof:}
First assume that $R$ is algebraically closed of characteristic $0$.
Let $\L=\bigoplus_{i=0}^k \L_i$ be a direct sum of simple 
Lie superalgebras as in the proposition
with Casimir element $\omega=\oplus_i \omega_i$.
The following map is a well-defined morphism of abelian groups:

$$
\pi_k:G_0(\L_1)\otimes_{\Z}\ldots\otimes_{\Z} G_0(\L_k)\longrightarrow G_0(\L)\ ,\ \
[\rho_1]\otimes_\Z\ldots\otimes_\Z[\rho_k]\mapsto[\rho_1\otimes\ldots\otimes
\rho_k].
$$

In Theorem~8 of \cite{Kac} simple $\L_i$-modules are classified by
highest weights.
This classification implies a classification of the simple $\L$-modules by highest 
weights. 
If we choose generators for each $\L_i$ as in Theorem~8 of \cite{Kac} and
$v_i$ are highest weight vectors with respect to the chosen generators for $\L_i$ 
with weight $\alpha_i$, then $v_1\otimes\ldots\otimes v_k$ is a highest weight vector
for $\L$ with weight $\oplus_i \alpha_i$ because Cartan subalgebras are by
definition of degree $\overline{0}$.
Using a suitable order on the set of weights we can prove by induction that
$\pi_k$ is surjective.
\footnote{The theorem of Jordan-H\"{o}lder implies 
that $\pi_k$ is an isomorphism but we will not use this.}

Let $m_k:{(A^*)}^{\otimes k}\longrightarrow A^*$ denote the $k$-fold product of weight
systems.
Exer\-cise~6.33 of \cite{BN1} holds also true for Lie superalgebras. Using this we see
that the following diagram commutes:

\setlength{\unitlength}{1pt}
$$
\begin{picture}(150,80)(-20,-25)
\put(-37,50){\makebox(0,0){$G_0(\L_1)\otimes_\Z\ldots\otimes_\Z G_0(\L_k)$}}
\put(158,50){\makebox(0,0){$G_0(\L)$}}
\put(75,-19){\makebox(0,0){$A^*$}}
\put(3,41){\vector(1,-1){57}}
\put(148,41){\vector(-1,-1){57}}
\put(30,49){\vector(1,0){107}}
\put(75,57){\makebox(0,0){$\pi_k$}}
\put(-130,11){\makebox{$m_k\circ(W_{\L_1,\omega_1,\cdot}\otimes\ldots\otimes
W_{\L_k,\omega_k,\cdot})$}}
\put(128,11){\makebox{$W_{\L,\omega,\cdot}$}}
\end{picture}\nopagebreak
$$

Now let $\rho$ be some finite-dimensional representation of $\L$. We choose an element
$\sum_i n_i [\rho_{i,1}]\otimes\ldots\otimes[\rho_{i,k}]\in \pi_k^{-1}([\rho])$.
The elements $P_d\in A_d$ from Proposition~\ref{simpleLie} satisfy 
$W_{\L_j,\omega_j,\rho_{i,j}}(P_d)=0$
for all $i,j$. This implies

$$
W_{\L,\omega,\rho}(P_d)=\sum_i n_i \left(
W_{\L_1,\omega_1,\rho_{i,1}}\cdot\ldots\cdot W_{\L_k,\omega_k,\rho_{i,k}}\right)
(P_d)=0,
$$

because $P_d$ is primitive.

The structure of $G_0(\L)$ may change when the field $R$ is not algebraically closed,
but a simple argument will allow us to reduce this case to the
algebraically closed one (compare also Exercise~6.32 of \cite{BN1}).
Let $R$ be arbitrary of characteristic $0$ and let $\overline{R}$ be the
algebraic closure of $R$. The Lie algebra $\L=\bigoplus_i \L_i$ with Casimir 
$\omega=\oplus_i \omega_i$ and $A$ are now defined over $R$ and 
$\overline{R}\otimes \L$ ($\overline{R}\otimes A$) denotes the extension of 
scalars from $R$ to $\overline{R}$.
Then the following diagram commutes for all $R$-linear projections 
$p:\overline{R}\longrightarrow R$.

\setlength{\unitlength}{20pt}
$$
\begin{picture}(12,6)(-6,-3)
\put(-5,2){\makebox(0,0){$G_0(\L)$}}
\put(5,2){\makebox(0,0){$G_0(\overline{R}\otimes\L)$}}
\put(-4,2){\vector(1,0){7.5}}
\put(0,2.5){\makebox(0,0){$[\rho]\mapsto[id_{\overline{R}}\otimes \rho]$}}
\put(-5,-2){\makebox(0,0){$\Hom_R(A,R)$}}
\put(5,-2){\makebox(0,0){$\Hom_{\overline{R}}(\overline{R}\otimes A,\overline{R})$}}
\put(2.5,-2){\vector(-1,0){5.5}}
\put(0,-2.5){\makebox(0,0){$p\circ (\varphi\vert A)\leftarrow\hspace{-6pt}\raisebox{1pt}{$\scriptscriptstyle |$}\;\varphi$}}
\put(-5,1.5){\vector(0,-1){3}}
\put(-6,0){\makebox(0,0){$W_{\L,\omega,\cdot}$}}
\put(5,1.5){\vector(0,-1){3}}
\put(6.5,0){\makebox(0,0){$W_{\overline{R}\otimes\L,\omega,\cdot}$}}
\put(0,0){\makebox(0,0){$\overline{R}\otimes\Hom_R(A,R)$}}
\put(-4.5,1.5){\vector(3,-1){3}}
\put(-1,1){\makebox(0,0){$id_{\overline{R}}\otimes W_{\L,\omega,\cdot}$}}
\put(0.5,-0.5){\vector(3,-1){3}}
\put(2.9,-1){\makebox(0,0){$\cong$}}
\end{picture}\nopagebreak
$$

This completes the proof because the element $P_d$ of the proof of 
Proposition~\ref{simpleLie} is defined over $\Q$.
$\Box$

\bigskip

{\small \tt
Institut de Recherche Math\'{e}matique Avanc\'{e}e

Universit\'{e} Louis Pasteur - C.N.R.S.

7 rue Ren\'{e} Descartes

67084 Strasbourg Cedex, France

Email: lieberum@math.u-strasbg.fr
}
}
\end{document}